\documentclass[12pt]{amsart}
\usepackage{amsfonts}
\usepackage{latexsym,amsmath,amsthm,amssymb,graphicx,rawfonts}

\numberwithin{equation}{section}
\allowdisplaybreaks

\def\endproof{$\hfill\Box$\\}
\def\s{\,\,\,\,}
\def\R{\mathbb{R}}
\def\C{\mathbb{C}}

\def\N{\mathbb{N}}
\def\mv{2.0ex}

\numberwithin{equation}{section}
\newtheorem{theorem}{Theorem}[section]
\newtheorem{defi}{Definition}
\newtheorem{Thm}{Theorem}
\newtheorem{lem}[theorem]{Lemma}
\newtheorem{thm}[theorem]{Theorem}
\newtheorem{pro}[theorem]{Proposition}
\newtheorem{cor}[theorem]{Corollary}
\newtheorem{rem}[theorem]{Remark}

\def\dint{\displaystyle{\int}}
\def\dis{\displaystyle}
\def\lan{\langle}
\def\ran{\rangle}
\def\diam{{\rm diam\,}}

\begin{document}
\title[Bubble tree, Conformal mapping, Willmore functional]
{\bf Bubble tree of a class of conformal mappings and applications to Willmore functional}
\author{Jingyi Chen \& Yuxiang Li}
\date{}
\thanks{The first author is partially supported by NSERC. The second author was partially supported by NSERC during his visit to UBC in the spring of 2011 where most part of this work was done}
\maketitle

\begin{abstract}
We develop a bubble tree construction and prove compactness results for $W^{2,2}$ branched conformal immersions of closed Riemann surfaces, with varying conformal structures whose limit may degenerate, in a compact Riemannian manifold with uniformly bounded areas and Willmore energies. The compactness property is applied to construct Willmore type surfaces in compact Riemannian manifolds. This includes (a) existence of a Willmore 2-sphere in ${\mathbb S}^n$ with at least 2 nonremovable singular points (b) existence of minimizers of the Willmore functional with prescribed area in a compact manifold $N$ provided (i) the area is small when genus is 0 and (ii) the area is close to that of the area minimizing surface of Schoen-Yau and Sacks-Uhlenbeck in the homotopy class of an incompressible map from a surface of positive genus to $N$ and
$\pi_2(N)$ is trivial (c) existence of smooth minimizers of the Willmore functional if a Douglas type condition is satisfied.
\end{abstract}

\section{Introduction}
Let $\Sigma$ be a smooth Riemann surface and $\ f : \Sigma\rightarrow \R^n\ $ be a smooth immersion.
The Willmore functional of $f$ is defined by
\begin{displaymath}
    W(f) = \frac{1}{4} \int_\Sigma |H_f|^2 d \mu_{f}
\end{displaymath}
where $H_f=\Delta_{g_f}f$ denotes the mean curvature vector of $f$,
 and $\Delta_{g_f}$ is the Laplace operator in the induced metric $g_f$ and $d\mu_f$ the induced area element on $\Sigma$.

For a sequence of immersions $f_k$ of a compact surfaces $\Sigma$ in ${\mathbb R}^n$ with uniformly bounded areas $\mu(f_k)$ and Willmore functionals $W(f_k)$, a subsequence of $\Sigma_k = f_k(\Sigma)$
converges, as Radon measures, to a two dimensional integral varifold, by Allard's integral compactness theorem. The second fundamental forms $A_{f_k}$ are uniformly bounded in the $L^2$-norm as
$$
\int_\Sigma|A_{f_k}|^2d\mu_{f_k}=4W(f_k)-4\pi\chi(\Sigma)
$$
from the Gauss equation and the Gauss-Bonnet formula. In general $\|f_k\|_{W^{2,2}}$ are not uniformly bounded: we can find diffeomorphisms $\phi_k$ from $\Sigma$ to $\Sigma$ such that $f_k=f\circ \phi_k$ diverge in $C^0$, while a uniform bound on $\|f_k\|_{W^{2,2}}$ would imply sequential convergence in $C^0$ (in fact $C^\alpha,0<\alpha<1$) norm by the Rellich-Kondrachov embedding theorem.

A recent advance in understanding the limit process is given in \cite{K-L}, where
 each $f_k$  is a conformal immersion
from a Riemann surface $(\Sigma,h_k)$ into $\R^n$ and $h_k$ is the
smooth metric of constant curvature:
\begin{equation}\label{metric-c}
\begin{array}{l}
h_k\mbox{ has Gauss curvature }\pm 1, \mbox{ or }
(\Sigma,h_k)=\C/\{1, a+bi\}
\mbox{ with }\\
-\frac{1}{2}
<a\leq\frac{1}{2},\s b\geq 0,
a^2+b^2\geq 1
\mbox{ and }a\geq 0
\mbox{ whenever }a^2+b^2=1.
\end{array}
\end{equation}
There are two  reasons to use conformal immersions. One is  that the conformal
diffeomorphism group of $(\Sigma,h_k)$ is rather small comparing with the group of diffeomorphisms.
Secondly, if we set $g_{f_k}=e^{2u_k} g_{euc}$ in an
isothermal coordinate system, then we can estimate
$\|u_k\|_{L^\infty}$ from the compensated compactness property  of
$K_{f_k}e^{2u_k}$. Thus it is possible to get an upper bound of
$\|f_k\|_{W^{2,2}}$ via the equation $\Delta_{h_k}f_k=H_{f_k}$.
When the conformal structures determined by $f_k$ do not go to the boundary of the moduli space, convergence of $f_k$ is treated in \cite{K-L}:  if the  conformal classes induced by $f_k$
converge in the moduli space, then there exist Mobius transformations
 $\sigma_k$, such that $\sigma_k\circ f_k$ converge locally  in
weak $W^{2,2}$ sense on $\Sigma$ minus finitely many concentration points.
The weak limit $f_0$ is a $W^{2,2}$ branched conformal immersion.

The $W^{2,2}$ conformal immersions and  $W^{2,2}$ branched conformal
immersions are as follows:
\begin{defi}\label{defconformalimmersion}{\em
Let $(\Sigma,h)$ be a connected Riemann surface with $h$ satisfies (\ref{metric-c}). A map $f\in W^{2,2}(\Sigma,h,\mathbb{R}^n)$
is called a {\em conformal immersion of $(\Sigma,h)$}, if
$$d f\otimes d f=e^{2u}h\s \hbox{with}\s \|u\|_{L^\infty(\Sigma)}<+\infty.$$
We denote the set of all such immersions by $W^{2,2}_{conf}(\Sigma,h,\R^n)$.
 It can be shown that for
$f\in W^{2,2}_{conf}(\Sigma,\R^n)$ the corresponding $u$ is continuous (see Appendix). When $f\in W^{2,2}_{loc}(\Sigma,h,\R^n)$  with
$d f\otimes d f=e^{2u}h$ and $u\in L^\infty_{
loc}(\Sigma)$, we say $f\in W^{2,2}_{conf,loc}(\Sigma,h,\R^n)$.}
\end{defi}

\begin{defi}{\em
We say $f$ is a {\em $W^{2,2}$
branched conformal  immersion of $(\Sigma,h)$ with possible branch points
$x_1, \dots ,x_m$}, if $f\in W^{2,2}_{conf,loc}(\Sigma\backslash
\{x_1,\dots,x_m\},h,\R^n)$ and
$$\int_{\Sigma\backslash\{x_1, \dots ,x_m\}}(1+|A_{f}|^2)
d\mu_f<+\infty.$$
The set of $W^{2,2}$ branched conformal immersions is denoted by $W^{2,2}_{b,c}(\Sigma,h,\R^n)$.
We say $f\in \widetilde{W}^{2,2}(\Sigma,\R^n)$, if
we can find a smooth metric $h$ satisfying (\ref{metric-c}) over $\Sigma$, such that
$f\in W^{2,2}_{b,c}(\Sigma,h,\R^n)$. }
\end{defi}

The first part of the paper is a study of a sequence of $W^{2,2}$ branched conformal immersions and the main goal is to establish compactness in Hausdorff distance for such immersions with uniformly bounded areas and Willmore functionals (cf. Theorem \ref{main}).


Our compactness result holds not only when $h_k$ converges, but also when
the conformal classes $c_k$ of $h_k$ diverge in the moduli space $\mathcal{M}_g$.
Bubbles develop near points where the Willmore energy concentrates, and if $c_k$ go to a point in
the boundary $\overline{\mathcal M}_g\backslash\mathcal{M}_g$ additional complication arises as the topology of the limit may be different from that of $\Sigma$ and stratified surfaces are used as possible limits.
The main idea to deal with degenerating conformal structures in the limit process is as follows.
First, pull the immersions $f_k$ to the immersions from components of
$\Sigma_0$, by composing $f_k$ with diffeomorphisms from the regular parts of the limit $\Sigma_0$ of $(\Sigma,h_k)$ in $\overline{\mathcal{M}_p}$. Then we study
the limit of $f_k$ and construct bubble trees at the energy
concentration points and collars, and investigate behavior between bubbles. In particular, we will prove that there is no loss of measure in the limit and there is no neck between the bubbles.  Then the limit $f_0$ of $f_k$ is a union of conformal maps from some components $\overline{\Sigma_0^1}, \dots ,
\overline{\Sigma_0^m}$ of $\Sigma_0$ (we delete those components
whose images are points) and finitely many 2-spheres $S_1,\dots,S_l$ into $\R^n$.
``no neck" means that we can glue  $\overline{\Sigma_0^i}$'s
and $S_j$'s to a stratified surface $\Sigma_\infty$ (see definition below), and
$f_0$ is a continuous map from $\Sigma_\infty$ into
$\R^n$. Then we will apply a result of H\'elein \cite{He} and a
removable singularity theorem in \cite{K-L} to
 show that for a sequence of branched conformal immersions with uniformly bounded measures and Willmore functionals, the limit we get in
 section 2 is in fact a branched conformal immersion
 of a stratified surface.

We point out that  the ``no loss of measure'' and ``no neck'' phenomenon must occur whenever the following two equations hold:
\begin{equation}\label{intro.mean}
-\Delta f_k=\frac{1}{2}|\nabla f_k|^2H_k,\s \mbox{with}\s \sup_k\int|\nabla f_k|^2(1+|H_k|^2)<\infty,
\end{equation}
\begin{equation}\label{conformal}
\partial f_k\otimes \partial f_k=0\,\,\,(\hbox{weakly conformal})
\end{equation}
where $\Delta,\nabla,\partial =\partial/\partial z$ are the operators in $h_k$ and its conformal structure $c_k$.
In section 2, we study the blowup behavior of a sequence
satisfies \eqref{intro.mean} and \eqref{conformal}.

The equation \eqref{intro.mean} looks similar to the equation of harmonic maps
$$-\Delta u=A(u)(du,du).$$
In fact, the arguments in section 2 are originated from the
``energy identity'' and ``no neck'' arguments of
harmonic maps \cite{C-T,  D-T, La, L-W, Lin-W, P, Q, Q-T, S-U1}.
When conformal structures go to the boundary of ${\mathcal M}_g$, non-trivial necks exist for harmonic map ( \cite{C-T, C-L-W, P, Z}); in our case, however, there is no non-trivial neck due to conformality although
 \eqref{intro.mean}  is much weaker than the harmonic map equation.

\begin{defi}{\em
Let $(\Sigma, d)$ be a connected compact metric space. We say $\Sigma$ is
a  {\em stratified surface with singular set $P$}
if $P\subset\Sigma$ is finite set such that

1. $(\Sigma\backslash P,d)$ is a smooth Riemann
surface without boundary (possibly disconnected) and $d$ is a smooth
metric $h=d|_{\Sigma\backslash P}$, and

2. For each $p\in P$, there is $\delta$, such that
$B_\delta(p)\cap P=\{p\}$ and $B_\delta(p)\backslash\{p\}
=\bigcup\limits_{i=1}^{m(p)}\Omega_i$,
where $1<m(p)<+\infty$, and each $\Omega_i$ is topologically
a disk with its center deleted. Moreover, on each $\Omega_i$,
$h$ can be extended to be a smooth metric on the disk.
}
\end{defi}

The genus of $\Sigma$ is defined by
$$g(\Sigma)=\frac{2-\chi(\Sigma)+\sum\limits_{p\in P}(m(p)-1)}{2}.$$
When $g(\Sigma)=0$, $\Sigma$ is called a stratified sphere.  A stratified
surface with singular set $P=\emptyset$ is a smooth Riemann surface.

\begin{center}
\includegraphics{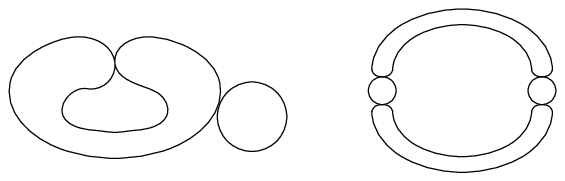}

Figure 1. Stratified torus
\end{center}

For a stratified surface $\Sigma$ with singular set $P$, we can write
$\Sigma\backslash P=\bigcup_i\Sigma^i$
where $\Sigma^i$'s are the disjoint connected components of $\Sigma$, and each
$\Sigma^i$ is a punctured Riemann surface when there are more than one components. The
topological closure of $\Sigma^i$ is denoted by $\Sigma_i$, so as a point-set $\Sigma=\bigcup_i \Sigma_i$.
By 2 in the above definition, each component $\Sigma^i$  can be extended
to a closed Riemann surface $\overline{\Sigma^i}$ by adding finitely many points. To illustrate the difference of these notations, take, for example, the stratified torus on the left in Figure 1: $P$ contains two points,
$\Sigma^1$ is the ``torus" with two points deleted and $\Sigma^2$ is a 2-sphere with one point removed,
$\Sigma_1$ is the ``torus" and $\Sigma_2$ is the 2-sphere, while $\overline{\Sigma^1}$ is a Riemann sphere (adding 3 points at the punctures) and $\overline{\Sigma^2}$ is also a Riemann sphere (adding 1 point at the puncture).

When $\Sigma$ is a stratified surface we define $f\in
W^{2,2}_{b,c}(\Sigma,\R^n)$ if $f$ is a $W^{2,2}$
non-trivial branched conformal immersion on each
$\overline{\Sigma^i}$.\\

We now state the main result in the first part of the paper:

\begin{Thm}\label{main}
Suppose that $\{f_k\}$ is a sequence of $W^{2,2}$  branched conformal
immersions of $(\Sigma,h_k)$ in $\R^n$, where $h_k$
satisfies \eqref{metric-c}. If
$f_k(\Sigma)\cap B_{R_0}\neq\emptyset$  for some fixed $R_0$ and
$$
\sup_k\left(\mu(f_k)+W(f_k)\right)<+\infty
$$
then  either $f_k$ converges to a point, or
there is a stratified surface $\Sigma_\infty$ with
$g(\Sigma_\infty)\leq g(\Sigma)$, a map
$f_0\in W^{2,2}_{b,c}(\Sigma_\infty,\R^n)$,
such that a subsequence of
$f_k(\Sigma)$ converges to  $f_0(\Sigma_\infty)$
in Hausdorff distance with
$$\mu(f_0)=\lim_{k\rightarrow+\infty}\mu(f_k)\s and\s W(f_0)\leq\lim_{k\rightarrow+\infty} W(f_k).$$
Moreover, if
$y_1,\dots,y_m\in f_k(\Sigma)$ for all $k$,
then
$y_1,\dots,y_m\in f_0(\Sigma_\infty).$
\end{Thm}

In fact, we will prove that $f_k$ converges to $f_0$ in the sense of bubble tree.
For each $k$, we can find a domain  $U_k$ of $\Sigma$ and
a domain $V_k$ of $\Sigma_\infty$, such that

1) $V_k\subset V_{k+1}$, and $P=\Sigma_\infty\backslash
\bigcup_kV_k$ is a finite set which contains all
singular points of $\Sigma_\infty$. Moreover,
$\Sigma_\infty\backslash V_k$ is a union of topological
disks and $H^1_1(\Sigma_\infty\backslash V_k)\rightarrow 0$,
where $H^1_1$ is the Hausdorff measure:
$$
H_1^1(S)=\inf\left\{\sum_{i=1}^\infty\mbox{\rm diam}(U_i):S\subset\bigcup_{i=1}^\infty
U_i,\s \mbox{\rm diam}(U_i)<1\right\}.
$$

2) $\Sigma\backslash U_k$ is a smooth surface with boundary, possibly disconnected, $H^1_1(f_k(\Sigma\backslash U_k))\rightarrow 0$. Moreover, $f_k(\Sigma\backslash U_k)$ converge to $P$ in Hausdorff distance.

3) There is a sequence of diffeomorphisms $\phi_k: V_k
\rightarrow U_k$, such that for any $\Omega\subset\subset
\Sigma_\infty\backslash P$,
$f_k\circ\phi_k$ converge in $W^{2,2}(\Omega,\R^n)$ weakly.

In Theorem \ref{main}, the singular points of $f_0(\Sigma_{\infty})$ arise in three ways: (a) the limit point to which a sequence of closed geodesics that are not null-homotopic in $f_k(\Sigma)$ pinches, (b) a bubble point of $f_k$, so belonging to a 2-sphere (the bubble), (c) a point where both (a) and (b) happen.



\vspace{0.3cm}

In the second part of the paper, we apply Theorem \ref{main} to obtain several existence results
of Willmore surfaces in compact Riemannian manifolds. Here we note that Theorem \ref{main} is applicable for surfaces immersed in a compact Riemannian manifold $N$. To see this,  for $\Sigma$ immersed in $N$ which is isometrically embedded in $ \R^n$,  direct calculation shows that
the Willmore functional of $\Sigma$ in $\R^n$ is dominated by its Willmore functional in $N$ together with the area $\mu(\Sigma)$, see Lemma \ref{H}.

We first consider 2-spheres immersed in the round unit sphere ${\mathbb S}^n, n\geq 3$. Fix at least two distinct points $y_1, \dots , y_m, m\geq 2$ on ${\mathbb S}^n$. Define
$$
\beta^n_0(y_1, \dots, y_m) = \inf \left\{W_n(f): f\in W^{2,2}_{conf}(S^2,{\mathbb S}^n), y_1,\dots,y_m\in f({ S}^2) \right\}
$$
where $W_n(f)=\int_{S^2}\left(1+ \frac{1}{4}\left|H_f\right|^2\right)d\mu_f$ and $H_f$ is the mean curvature vector of $f(S^2)$ in ${\mathbb S}^n$.

We show
\begin{Thm}\label{sphere}
If $\beta^n_0(y_1,\dots,y_m)<8\pi$, then there is a $W^{2,2}$ conformal immersion of $S^2$ in ${\mathbb S}^n$
without self-intersections realizing $\beta^n_0(y_1,\dots,y_m)$. For any $\epsilon > 0$, there exists a Willmore sphere in ${\mathbb S}^n$ with $  W _n(f ) < 4\pi + \epsilon$, which
has at least 2 nonremovable singularities.
\end{Thm}

By results in \cite{K-S2},  \cite{R}, a singular point of a
Willmore surface with density $\theta^2<2$ in $\R^n$
can be removed if its residue is 0.
Kuwert and Sch\"atzle also point out that the removability can
not be true generally, for example, 0 is the true
singular point of an inverted half catenoid (\cite{K-S2}, P. 337).
The second statement in Theorem \ref{sphere} provides examples of embedded Willmore surface which
has a nonremovable singular point with density $\theta^2=1$, and it is an application of the first statement with five points prescribed in ${\mathbb S}^n$.

We then consider minimizers of the Willmore functional subject to area constraint. A fundamental existence result for incompressible minimal surfaces due to Schoen-Yau \cite{S-Y} and Sacks-Uhlenbeck \cite{S-U2}
asserts:
If $\varphi$ induces an injection from the fundamental groups to $\Sigma$ and $N$, then there is a branched minimal immersion $f: \Sigma\rightarrow N$ so that $f$ induces the same map between fundamental groups as $\varphi$ and $f$ has least area among all such maps. We denote the area of the minimizer by $a_\varphi$.

\begin{Thm}\label{area}
Let $N$ be a compact Riemannian manifold and $\Sigma$ be a closed surface of genus $g$. Then

(1) For $\beta_0(N,a) = \inf\{W(f): \mu(f) = a>0, f\in W^{2,2}_{conf}(S^2,N)\}$, $\lim_{a\to 0}\beta_0(N,a)=4\pi$, and there is an embedding realizing $\beta_0(N,a)$ for all sufficiently small $a$.

(2) Suppose $\varphi:\Sigma\to N$ induces an injection $\varphi_{\#}:\pi_1(\Sigma)\to\pi_1(N)$ and $\pi_2(N)=0$. Let $\beta_g(N,a,\varphi)=\inf\{W(f):f\in\widetilde{W}^{2,2}(\Sigma,N),\mu(f)=a,  f \hbox{ is homotopic to}\,\,\, \varphi\}$.
Then there is $\delta > 0$, such that for any $a \in [a_\varphi , a_\varphi + \delta)$, there is a
branched conformal immersion $f$ of $(\Sigma,h)$ attaining $\beta_g(N, a, \varphi)$.
Moreover, when $\dim N = 3$,  $f$ is an immersion for small $\delta$.
\end{Thm}

For $\beta_0(N,a)$,  Lamm and Metzger showed in \cite{La-M}
that if it is attained by a surface with
positive mean curvature in the sufficiently small
geodesic ball around a point $p$, then the scalar
curvature of $N$ must have a critical point at $p$.

When $N$ has negative sectional curvature, the area of an immersed surface is dominated by the Willmore functional. We now describe a sufficient condition of Douglas type for existence.
Let $S(g)$ be the set of connected stratified Riemann surfaces $\Sigma=\bigcup_i\Sigma_i$ satisfying (a) genus of $\Sigma_i < g$ if $g > 0$ and (b) $i > 1$ if $g = 0$. Note that a surface in $S(g)$ has genus at most $g$ and
smooth surfaces of genus $g$ are not in $S(g)$. Isometrically embed $N$ into ${\mathbb R}^n$.  Define
$$
\alpha^*(g) = \inf \{ W(f) : f  \in W^{2,2}_{b,c} ( \Sigma,
\mathbb{R}^n), \Sigma \in {S}(g) \}
$$
$$
\alpha(g) = \inf \{ W(f) : f \in W^{2,2}_{b,c}(\Sigma,{\mathbb R}^n), \hbox{$\Sigma$ is
a smooth surface of genus $g$}\}.
$$
Similarly, for $0<a<\infty$, define
\begin{eqnarray*}
\gamma^*(g,a)&=&\inf \{ W(f,\Sigma,\R^n): f\in W^{2,2}_{b,c}(\Sigma,\R^n), f(\Sigma)\subset N,\Sigma\in S(g), \mu(f(\Sigma))\leq a\} \\
\gamma(g,a)&=& \inf \{W(f,\Sigma,\R^n): f\in W^{2,2}_{b,c}(\Sigma,\R^n), f(\Sigma)\subset N,\Sigma\in{\mathcal M}_g,
\mu(f(\Sigma))\leq a\}.
\end{eqnarray*}
\begin{Thm}
Let $N$ be a compact Riemannian manifold. If $0<\alpha(g)<\alpha^*(g)$ and $N$ has negative sectional curvature, then there is a $W^{2,2}$ branched conformal immersion $f$ from a smooth closed Riemann surface of genus $g$ with $W(f)=\alpha(g)$. If $0<\gamma(g,a)<\gamma^*(g,a)$ then there is a $W^{2,2}$ branched conformal immersion $f$ from a smooth closed Riemann surface of genus $g$ with $W(f)=\gamma(g)$.
\end{Thm}

\section{blowup analysis - energy identity and absence of neck}
Let $(\Sigma,h)$ be a smooth Riemann surface which may not be compact, where $h$
is the metric compatible with the complex structure
of $\Sigma$.  For given $p>1$ and $R>0$, let $\mathcal{F}^p(\Sigma,h, R)$ be the set of mappings $f:
\Sigma\rightarrow\R^n$ which satisfy
\begin{enumerate}
\item $f\in W^{2,p}_{loc}(\Sigma,h)$;
\item $f(\Sigma)$ is contained in the closed ball centered at the origin with radius $R$ in ${\mathbb R}^n$;
\item $\Delta_hf=F(f)$ with
$|F(f)|\leq \beta\, |\nabla_hf|^2$ a.e. on $\Sigma$, where $\beta$
is a nonnegative measurable function on $\Sigma$ with
$$
\int_\Sigma\beta^2|\nabla_hf|^2d\mu_h<+\infty.
$$
\end{enumerate}
When $f\in\mathcal{F}^p(\Sigma,h,R)$, we introduce a notation by
$$
H(f)=\left\{\begin{array}{ll}
                     2\frac{\Delta_h f}{|\nabla_h f|^2},& \hbox{if $|\nabla_h f|\neq 0$}\\
                     0,& \hbox{if $|\nabla_hf|=0$}.
              \end{array}\right.
              $$
The Willmore functional of $f$ is defined to be
$$
W(f)=\frac{1}{4}\int_\Sigma|H(f)|^2|\nabla_hf|^2d\mu_h.
$$
That  $W(f)<\infty$ for $f\in\mathcal{F}^p(\Sigma,h,R)$ follows from (3) as
$$
\Delta_h f=F(f)=\frac{1}{2}H(f)|\nabla_h f|^2.
$$
We denote by $\mathcal{F}^p_{conf}(\Sigma,h,R)$ the set of $f\in\mathcal{F}^p(\Sigma,h,R)$ and $f$ is weakly conformal a.e., i.e. $\partial f\otimes\partial f=0$ almost everywhere on $\Sigma$, where $\partial f=\frac{\partial f}{\partial z}dz$ in a local complex coordinate system on $\Sigma$.

Note that when $f$ is a smooth conformal immersion $H(f), \frac{1}{2}|\nabla_h f|^2d\mu_h, W(f)$ are the
mean curvature vector,  the area element  and the Willmore functional of $f(\Sigma)$, respectively.

By the Kondrachov embedding theorem, functions in ${\mathcal F}^p(\Sigma,h,R)$ are also in $W^{1,2}$.
The right hand side of the equation $\Delta_h f = F(f)$ is not necessarily in $L^2$ under the assumption (3).

We point out that $H(f)$, $\mathcal{F}^p,\mathcal{F}^p_{conf}$ are conformal invariant, in the sense that if $h'=e^{2u}h$ for some smooth function $u$ on $\Sigma$, we always have
$$
H_h(f)=H_{h'}(f), \s
\mathcal{F}^p(\Sigma,h,R)=\mathcal{F}^p(\Sigma,h',R),\s \mathcal{F}^p_{conf}(\Sigma,h,R)=\mathcal{F}^p_{conf}(\Sigma,h',R).
$$
Thus we may select preferred metrics $h$, e.g. the ones with constant curvature.

In this section, we will study regularity, compactness and the blowup behavior of
a sequence $\{f_k\}\subset\mathcal{F}^p$.

\subsection{$\epsilon$-regularity, removable singularity and weak limit}
In this subsection, we will show that some well-known
results for harmonic maps still hold for mappings in $\mathcal{F}^p$.

Let $D$ be the unit 2-disk centered at 0. For simplicity, write $\mathcal{F}^p(D,dx^2+dy^2,R)$
as $\mathcal{F}^p(D,R)$.
\begin{pro}\label{epsilon} ($\epsilon$-regularity) There is an $\epsilon_0$ such that for any $f\in \mathcal{F}^p(D,R)$, $1<p<2$, if $W(f)<\epsilon_0^2$, then
$$
\|\nabla f\|_{W^{1,p}(D_\frac{1}{2})}\leq C\,\|\nabla f\|_{L^2(D)}.
$$
\end{pro}

\proof Set $\bar{f}=\frac{1}{|D|}\int_{D}
fd\sigma$ and let $\eta$ be a cut-off function which is 1 in
$D_{1/2}$, 0 in $D\backslash D_{3/4}$ and $0\leq \eta\leq 1$. Then for the equation
$$
\Delta \left(\eta(f-\bar{f})\right)=(f-\bar{f})\Delta\eta+2\nabla\eta
\nabla f+\frac{1}{2}\eta H(f)|\nabla f|^2 := \phi
$$
we have
\begin{eqnarray*}
|\phi|&\leq& C_1\left(|f-\bar{f}|+|\nabla f|\right)+\frac{1}{2}\eta \left|H(f)\right||\nabla f|^2\\
&\leq&C_1\left(|f-\bar{f}|+|\nabla f|\right)+C_2\left|H(f)\right||\nabla f|\left(|\nabla\left(\eta(f-\bar{f})\right)|+|f-\bar{f}|\right)
\end{eqnarray*}
since
$$
\begin{array}{lll}
\frac{1}{2}\eta \left|H(f)\right||\nabla f|^2&=&\frac{1}{2}\eta \left|H(f)\right|\nabla(f-\bar{f})\nabla f\\[\mv]
&=&\frac{1}{2}\left|H(f)\right|\nabla\left(\eta(f-\bar{f})\right)\nabla f-\frac{1}{2}\left|H(f)\right|(f-\bar{f})\nabla\eta\nabla f\\[\mv]
&\leq& C_2\left|H(f)\right||\nabla f|\left(|\nabla\left(\eta(f-\bar{f})\right)|+|f-\bar{f}|\right).
\end{array}
$$
By the $L^p$ estimates for elliptic equations,
\begin{eqnarray*}
\left\|\eta(f-\bar{f})\right\|_{W^{2,p}(D)}
   &\leq& C_3\left(\left\|f-\bar{f}\right\|_{L^p(D)}+\left\|\nabla f\right\|_{L^p(D)}\right.\\
   &&\left.+\left\|H(f)|\nabla f|\left(|\nabla\left(\eta(f-\bar{f})\right)|+
    |f-\bar{f}|\right)\right\|_{L^p(D)}\right).
\end{eqnarray*}
For $1<p<2$, the H\"older inequality and the Sobolev inequality imply
$$
\begin{array}{lll}
&&\left\|H(f)|\nabla f|\left(\left|\nabla\left(\eta(f-\bar{f})\right)\right|+
    \left|f-\bar{f}\right|\right)\right\|_{L^p(D)}\\
&&\leq\|H(f)\nabla f\|_{L^2(D)}\left(\|\nabla\left(\eta(f-\bar{f})\right)\|_{L^\frac{2p}{2-p}(D)}
  +\|f-\bar{f}\|_{L^\frac{2p}{2-p}(D)}\right)\\[\mv]
&&\leq \epsilon_0\,C_4\,\|\eta(f-\bar{f})\|_{W^{2,p}(D)}+\epsilon_0\,C_5\,\|f-\bar{f}\|_{W^{1,p}(D)}
\end{array}
$$
since $W(f)<\epsilon_0$.
Applying the Poincar\'e inequality and noting $1<p<2$, we get
$$\|f-\bar{f}\|_{L^p(D)}+\|\nabla f\|_{L^p(D)}+ \epsilon_0\,C_5\,\|f-\bar{f}\|_{W^{1,p}(D)}
\leq C_6\,\|\nabla f\|_{L^2(D)}.$$
Choose $\epsilon_0$ so that  $C_3C_4\,\epsilon_0<1/2$, then we  get
$$\|\eta(f-\bar{f})\|_{W^{2,p}(D)}<C_7\,\|\nabla f\|_{L^2(D)}
$$
which completes the proof. \endproof

\begin{pro}\label{gap} (Gap constant) Let $\Sigma$ be a closed surface.
There is a constant ${\epsilon}_1$ which depends
on $\Sigma$ and $R$, such that for any
$f\in\mathcal{F}^p(\Sigma,h,R)$ where $1<p<2$, if $W(f)<{\epsilon}_1^2$, then $f$ is constant.
\end{pro}

\proof Let $\bar{f}=\frac{1}{|\Sigma|}\int_\Sigma f$.
It follows from the equation
$$
\Delta_h(f-\bar{f})=\frac{1}{2}H(f)|\nabla f|^2
$$
that
$$\begin{array}{lll}
\dis\int_\Sigma|\nabla(f-\bar{f})|^2
    &\leq&\dis \frac{1}{2}\int_\Sigma|f-\bar{f}|\,|{H(f)}|\,|\nabla f|^2\\[\mv]
    &\leq&\dis \left(\int_\Sigma{H(f)}^2|\nabla f|^2\right)^\frac{1}{2}
     \left(\int_\Sigma|f-\bar{f}|^\frac{2p}{2p-2}\right)^\frac{2p-2}{2p}
     \left(\int_\Sigma|\nabla f|^\frac{2p}{2-p}\right)^\frac{2-p}{2p}\\[\mv]
     &\leq&\dis C_1W(f)^{\frac{1}{2}}\|\nabla f\|_{L^2(\Sigma)}\|\nabla f\|_{L^\frac{2p}{2-p}(\Sigma)}
\end{array}$$
where we used the Sobolev inequality, the Poincar\'e inequality and $1<p<2$. Then we get
$$
\|\nabla f\|_{L^2(\Sigma)}\leq C_1W(f)^{\frac{1}{2}}\|\nabla f\|_{L^\frac{2p}{2-p}(\Sigma)}.
$$
Using the Poincar\'e inequality and $1<p<2$ again, we have
$$\|f-\bar{f}\|_{L^p(\Sigma)}\leq C_2\|\nabla f\|_{L^2(
\Sigma)}\leq C_1C_2\,W(f)^{\frac{1}{2}}\|\nabla f\|_{L^\frac{2p}{2-p}(\Sigma)}.$$
Since
$$
\left\|\frac{1}{2}H(f)|\nabla f|^2\right\|_{L^p(\Sigma)}\leq \left(\int_\Sigma\frac{1}{4}{H(f)}^2|\nabla f|^2\right)^\frac{1}{2}
\left(\int_\Sigma|\nabla f|^\frac{2p}{2-p}\right)^\frac{2-p}
{2p}=W(f)^{\frac{1}{2}}\|\nabla f\|_{L^\frac{2p}{2-p}(\Sigma)},
$$
it follows from the $L^p$ estimates for elliptic equations that
\begin{eqnarray*}
\left\|f-\bar{f}\right\|_{W^{2,p}(\Sigma)}&\leq& C_3\,\left(\left\|\frac{1}{2}H(f)|\nabla f|^2\right\|_{L^p(\Sigma)}
+\left\|f-\bar{f}\right\|_{L^p(\Sigma)}\right)\\
&\leq& C_3(1+C_1C_2)W(f)^{\frac{1}{2}}\left\|\nabla f\right\|_{L^\frac{2p}{2-p}(\Sigma)}\\
&\leq&C_4W(f)^{\frac{1}{2}}\left\|f-\bar{f}\right\|_{W^{2,p}(\Sigma)}
\end{eqnarray*}
where the Sobolev inequality was used in the last step. By choosing $\epsilon_1<1/C_4$ we immediately
have $f=\overline{f}$. \endproof

We now derive a key estimate for later applications.
Set
$$
 E(f,Q(t))=\int_{Q(t)}|\nabla f|^2,\,\,\hbox{where}\,\,Q(t)=S^1\times[-t,t]
$$
and  denote  $\mathcal{F}^p(Q(t),dt^2+d\theta^2,R)$
by $\mathcal{F}^p(Q(t),R)$. We will prove the following energy decay
estimate:

\begin{pro}\label{key}(Decay estimate)
Let $f\in \mathcal{F}^p_{conf}(Q(T),R)$ with $T\gg 1,1<p<2$. Then there is a
constant $\epsilon_2<\epsilon_0$, where $\epsilon_0$ is the constant in Proposition \ref{epsilon}, such that if
$$
\sup_{t\in[-T,T-1]}W(f,S^1\times [t,t+1])<\epsilon^2\leq\epsilon_2^{2}
$$
then
$$\int_{Q(t)}|\nabla f|^2<{C}E(f,Q(T))e^{-(1-C\epsilon) (T-t)}
$$
for some positive constant $C$ independent of $T$ and $f$.
\end{pro}

\proof
 Define
$$
f^*(t)=\frac{1}{2\pi}\int_0^{2\pi}f(t,\theta)d\theta.
$$
We have
\begin{eqnarray}\label{f-star}
  \int_{Q(t)}\left|\frac{\partial f^*} {\partial t}\right|^2&=&
    \int_{-t}^{t}\int_0^{2\pi}\left(\frac{1}{2\pi}
     \int_0^{2\pi}\frac{\partial f}{\partial t}d\theta\right)^2
     d\theta dt\nonumber\\
         &\leq&\displaystyle\frac{1}{2\pi}\int_{-t}^{t}
      \left(\int_0^{2\pi}\left|\frac{\partial f}{\partial t}\right|^2d\theta
       \int_0^{2\pi}d\theta\right) dt\\
    &=&\int_{-t}^{t}\int_0^{2\pi}\left|\frac{\partial f}{\partial t}\right|^2dtd\theta\nonumber\\
    &=&\int_{Q(t)} \left|\frac{\partial f}{\partial t}\right|^2dtd\theta. \nonumber
      \end{eqnarray}
Then
\begin{equation}\label{id1}
\begin{array}{lll}\displaystyle\int_{Q(t)}\nabla (f-f^*)\nabla f&=&\displaystyle\int_{Q(t)}|\nabla f|^2-
      \displaystyle\int_{Q(t)}\frac{\partial f}{\partial t}\frac{\partial f^*}{\partial t}\\[\mv]
      &\geq&\displaystyle\int_{Q(t)}|\nabla f|^2-\frac{1}{2}\left(\displaystyle\int_{Q(t)}\left|\frac{\partial f}{\partial t}\right|^2
      +\displaystyle\int_{Q(t)}\left|\frac{\partial f^*}{\partial t}\right|^2\right)\\[\mv]
      &\geq&\displaystyle\int_{Q(t)}|\nabla f|^2-\int_{Q(t)}\left|\frac{\partial f}{\partial t}\right|^2\\[\mv]
      &=&\dis\frac{1}{2}\int_{Q(t)}|\nabla f|^2
  \end{array}
\end{equation}
where in the last step we used the fact that $|f_t|^2=|f_\theta|^2$ a.e. as $f$ is conformal a.e.
On the other hand,
\begin{equation}\label{id2}
 \begin{array}{lll}
    \displaystyle\int_{Q(t)}\nabla (f-f^*)\nabla f&=&-\displaystyle\int_{Q(t)}(f-f^*)\Delta f-
       \displaystyle\int_{\partial Q(t)}\frac{\partial f}{\partial t}(f-f^*)\\[\mv]
   &\leq&\displaystyle\int_{Q(t)}
    |f-f^*||\nabla f|^2\frac{|{H(f)}|}{2}
+\left|\displaystyle\int_{\partial Q(t)}
       \frac{\partial f}{\partial t}(f-f^*)\right|.
  \end{array}
\end{equation}

Let  $m\in [t,t+1)$ be an integer. Then for each $i=-m,
-m+1,\dots,m-1$, by (\ref{f-star}) and the hypothesis in the proposition
$$
 \sup_{t\in[-T,T-1]}\frac{1}{4}\int_{S^1\times[t,t+1]}|H(f)|^2<\epsilon^2\leq\epsilon_0^2
$$
it follows from Proposition \ref{epsilon} that
\begin{equation}\label{f^*}
\|f-f^*\|_{L^\infty(S^1\times[i,i+1])}\leq C\,\|\nabla f\|_{L^2(S^1\times[i-1,i+2])}.
\end{equation}
In fact, to see (\ref{f^*}), denote the average of $f$ over $S\times [i-1,i+2]$ by  $\overline{f}$
and observe that from Proposition \ref{epsilon}
$$
\| \nabla(f-\overline{f}) \|_{W^{1,p}(S^1\times [i,i+1])} = \| \nabla f \|_{W^{1,p}(S^1\times [i,i+1])}\leq C \|\nabla f\|_{L^2(S^1\times[i-1,i+2])}
$$
and from the Poincar\'e inequality
$$
\| f -\overline{f} \|_{L^2(S^1\times[i,i+1])} \leq \| f -\overline{f} \|_{L^2(S^1\times[i-1,i+2])} \leq C \| \nabla f \|_{L^2(S^1\times[i-1,i+2])}
$$
hence
$$
\| f-\overline{f} \|_{W^{2,p}(S^1\times[i,i+1])} \leq C \| \nabla f\|_{L^2(S^1\times[i-1,i+2])}.
$$
The Sobolev embedding theorem then implies
$$
\| f-\overline{f}\|_{C^{0}(S^1\times[i,i+1])} \leq C \| f-\overline{f} \|_{W^{2,p}(S^1\times[i,i+1])} \leq C\| \nabla f\|_{L^2(S^1\times[i-1,i+2])}.
$$
Therefore for any $t\in[i,i+1]$
$$
\left| f^*(t) - \overline{f} \right| = \left| \frac{1}{2\pi}\int^{2\pi}_0 \left(f(t,\theta)- \overline{f}\right) d\theta \right|
\leq C \| \nabla f\|_{L^2(S^1\times[i-1,i+2])}.
$$
It follows
$$
\|f-f^*\|_{C^0(S^1\times[i,i+1])}\leq
\|f-\bar{f}\|_{C^0(S^1\times[i,i+1])}+\|f^*-\bar{f}\|_{C^0(S^1\times[i,i+1])}
\leq C\|\nabla f\|_{L^2(S^1\times[i-1,i+2])}.
$$
Clearly, by the mean value theorem, $f-f^*$ equals 0 somewhere in $S^1\times[i,i+1]$, thus
$$
\| f-f^*\|_{L^\infty(S^1\times[i,i+1])} = \| f-f^*\|_{C^{0}(S^1\times[i,i+1])}\leq C\|\nabla f\|_{L^2(S^1\times[i-1,i+2])}
$$
which shows (\ref{f^*}) holds.

Then
$$\begin{array}{l}
\displaystyle\int_{S^1\times[i,i+1]}|f-f^*|
|\nabla f|^2\frac{|{H(f)}|}{2}\\[\mv]
   \begin{array}{lll}
     &\leq&\dis\left\|f-f^*\right\|_{L^\infty(S^1\times[i,i+1])}\times
     \left(W(f_k,S^1\times[i,i+1])  \int_{S^1\times[i,i+1]}|\nabla f|^2 \right)^{\frac{1}{2}}\\[\mv]
     &\leq&\dis C\epsilon\,\left(\int_{S^1\times[i-1,i+2]}|\nabla f|^2
     \int_{S^1\times[i,i+1]}|\nabla f|^2\right)^{\frac{1}{2}}\\[\mv]
     &\leq&\dis C\epsilon\int_{S^1\times[i-1,i+2]}|\nabla f|^2.
   \end{array}
\end{array}$$
Then
\begin{equation}\label{id4}
\begin{array}{lll}
\dint_{Q(t)}|f-f^*||\nabla f|^2{H(f)}
&\leq&\dis\sum\limits_{i=-m}^{m-1}\int_{S^1\times[i,i+1]}
|f-f^*||\nabla f|^2{H(f)}\\[\mv]
&\leq&\dis C\epsilon\sum\limits_{i=-m}^{m-1}
\int_{S^1\times
  [i-1,i+2]}|\nabla f|^2\\[\mv]
&\leq&\dis 3C\epsilon\int_{Q(t)}|\nabla f|^2\\[\mv]
&&\dis+C\epsilon\left(\int_{S^1\times[-m-1,-m]}|\nabla f|^2+
\int_{S^1\times[m,m+1]}|\nabla f|^2\right)\\[\mv]
&\leq&3C\epsilon\displaystyle\int_{Q(t+2)}|\nabla f|^2.
\end{array}
\end{equation}
From (\ref{id1}), (\ref{id2}), (\ref{id4}), we have
\begin{equation}\label{en}
\frac{1}{2}\int_{Q(t)} |\nabla f|^2 \leq \epsilon'\int_{Q(t+2)}|\nabla f|^2 +\left|\int_{\partial Q(t)}(f-f^*)\frac{\partial f}{\partial t}\right|
\end{equation}
where $\epsilon'=3C\epsilon/2$.
Moreover,
\begin{equation}\label{id3}
\begin{array}{lll}
 \dis\left|\int_{S^1\times\{t\}}\frac{\partial f}{\partial t}(f-f^*)\right|
         &\leq&\dis \left(\int_0^{2\pi}(f(\theta,t)-f^*(t))^2d\theta\right)^\frac{1}{2}
          \left(\int_0^{2\pi}\left|\frac{\partial f}{\partial t}(\theta,t)\right|^2d\theta\right)^\frac{1}{2}\\[\mv]
         &\leq&\dis\left(\int_0^{2\pi}\left|\frac{\partial f}{\partial \theta}(\theta,t)\right|^2d\theta\right)^\frac{1}{2}
           \left(\int_0^{2\pi}\left|\frac{\partial f}{\partial t}(\theta,t)\right|^2d\theta\right)^\frac{1}{2}\\[\mv]
         &=&\dis{\frac{1}{2}}\dint_{S^1\times\{t\}}|\nabla f|^2d\theta
  \end{array}
\end{equation}
here we used the Poincar\'e inequality on $S^1$ and
the fact that $|\frac{\partial f}{\partial t}|^2=|\frac{\partial f}{\partial \theta}|^2$ a.e. Let
$$
\varphi(t)=\frac{1}{2}\int_{Q(t)}|\nabla f|^2.
$$
By \eqref{en} and \eqref{id3}, we have
$$
\varphi(t)\leq \varphi'(t)+\epsilon'\varphi(t+2).
$$
Then
$$-(e^{-t}\varphi(t))'\leq\epsilon'\varphi(t+2)e^{-t},$$
and integrating the inequality from $t$ to $T-2$ leads to
\begin{equation}\label{id25}
\begin{array}{lll}
\dis e^{-t}\varphi(t)&\leq&\dis e^{-T+2}\varphi(T-2)+\epsilon'\int_t^{T-2}
\varphi(s+2)e^{-s}ds\\[\mv]
    &=&\dis e^{-T+2}\varphi(T-2)+\epsilon'\int_{t+2}^{T}\varphi(s)e^{-s+2}ds\\[\mv]
    &=&\dis e^{-T+2}\varphi(T-2)+\epsilon'e^2\int_{t+2}^{T-2}\varphi(s)e^{-s}ds
      +\epsilon'e^2\int_{T-2}^{T}\varphi(s)e^{-s}ds\\[\mv]
    &\leq&\dis e^{-T+2}\varphi(T)+\epsilon'e^2\int_t^{T-2}\varphi(s)e^{-s}ds+\epsilon'e^2\varphi(T)\left(e^{-T+2}-e^{-T}\right)
\end{array}
\end{equation}
as $\varphi(t)$ is increasing in $t$. Let $F(t)=\int_t^{T-2}\varphi(s)e^{-s}ds$ and
$\epsilon_2=\epsilon'e^2<1$. Now (\ref{id25}) leads to
$$
-F'(t)\leq 2\,\varphi(T)e^{-T+2}+\epsilon_2F(t)
$$
or equivalently
$$
\left(e^{\epsilon_2t}F(t)\right)'+2\varphi(T)e^{-T+2}e^{\epsilon_2t}\geq 0.
$$
Integrating over $[t,T-2]$ and noting $F(T-2)=0$, we have
\begin{equation}\label{ineq}
F(t)\leq \frac{2\varphi(T)}{\epsilon_2}e^{2-T}\left(e^{\epsilon_2(T-2)}-e^{\epsilon_2 t}\right)e^{-\epsilon_2t}.
\end{equation}
Substitute (\ref{ineq}) into  \eqref{id25}:
\begin{eqnarray*}
\varphi(t)&\leq& e^{2-T+t}\varphi(T)+2\varphi(T)e^{\epsilon_2(T-t-2)}e^{2-T+t}+\epsilon_2\varphi(T)e^{t-T+2}\\
&\leq& C\varphi(T)e^{(1-\epsilon_2)(T-t)}\\
&=& C\varphi(T)e^{(1-C\epsilon)(T-t)}
\end{eqnarray*}
for some positive constant $C$ independent of $T$ and $f$.
\endproof


\begin{pro}\label{rem}(Removability of point singularity)
Let $f\in \mathcal{F}^p_{conf}(D\backslash\{0\},R)$, where $1<p<2$. If $\int_D|\nabla f|^2<+\infty$,
then $f\in \mathcal{F}^{p'}_{conf}(D,R)$ for any $p'\in(1,\frac{4}{3})\cap(1,p]$.
\end{pro}

\proof We may assume that
$W(f)<\epsilon^2<\epsilon_2^{2}$,
otherwise, we can replace  $f$ with $f(\lambda x)$
for some $\lambda<1$.
Let $\phi:\R^1\times S^1\rightarrow\R^2$ be the conformal mapping
given by $r=e^{-t},\theta=\theta$.
Then $f'=f(\phi)$ is a map from $[0,+\infty)\times S^1$
into $\R^n$.
By translating $S^1\times [t-1,t+1]\subset S^1\times[0,2t]$ to $S^1\times [-1,1] \subset S^1\times[-t,t]$, from Proposition \ref{key} we conclude
$$
\int_{S^1\times[t-1,t+1]}|\nabla f'|^2\leq C_1e^{-\delta t},\s \hbox{where}\s \delta=1-C\epsilon.
$$
Then for any $r_k=e^{-k}$, we have $t_k=k$ and
\begin{equation}\label{in}
\int_{D_{r_{k-1}}\backslash D_{r_{k+1}}} |\nabla f|^2<C_1r_k^{-\delta }.
\end{equation}
Set $f_k(x)=f(r_kx)$. Applying Proposition \ref{epsilon} and (\ref{in}), we get
$$
\|\nabla f_k\|_{W^{1,p}(D_1\backslash D_{e^{-1}})}\leq C_2\,
\|\nabla f_k\|_{L^2(D_{e}\backslash D_{e^{-2}})}\leq
C_3\,r_k^\frac{\delta}{2}.
$$
By the Sobolev inequality, we have
$$
\left(\int_{D_1\backslash D_{e^{-1}}}|\nabla f_k|^q\right)^\frac{1}{q} \leq
C_4\,\|\nabla f_k\|_{W^{1,p}(D_1\backslash D_{e^{-1}})}\leq C_5\,r_k^{-\frac{\delta}{2}},\s where\s
q\leq \frac{2p}{2-p}.
$$
Then
$$\int_{D_{1}\backslash D_{e^{-1}}}|\nabla f_k|^q\leq C_6\,e^{-qk\frac{\delta}{2}}.$$
Since
$$
r_k^{2-q}\int_{D_1\backslash D_{e^{-1}}}|\nabla f_k|^q=\int_{D_{r_k}\backslash D_{r_{k+1}}}|\nabla f|^q,
$$
we have
$$
\int_{D_{r_k}\backslash D_{r_{k+1}}}|\nabla f|^q\leq C_6\,e^{-qk\frac{\delta}{2}+(q-2)k}=C_6\,e^{k(-2+q(1-\frac{\delta}{2}))}.
$$
When $q<4$, we can choose $\epsilon$ suitably
such that $q(1-\frac{\delta}{2})<2$, which yields
$$
\int_{D}|\nabla f|^q\leq C_6\sum_k2^{-qk\frac{\delta}{2}+(q-2)k}<C_7<\infty.
$$

For any $p'\in (1,\frac{4}{3})$, set $q=\frac{2p'}{2-p'}$, so $q\in (2,4)$. We have
$$
\int_D{H(f)}^{p'}|\nabla f|^{2p'}\leq\left(\int_D{H(f)}^2|\nabla f|^2\right)^\frac{p'}{2}
\left(\int_D|\nabla f|^q\right)^{\frac{p'}{q}}<C_8.
$$
Therefore, $F(f)\in L^{p'}(D)$ with $p'>1$ and then there exists $v$ which solves the equation
$$-\Delta v=F(f),\s v|_{\partial D}=0,$$
and $v\in W^{2,p'}(D)$. Obviously, $f-v$ is a harmonic function on $D\backslash\{0\}$ with
$$
\|\nabla(f-v)\|_{L^2(D)}+\|f-v\|_{L^2(D)}<+\infty.
$$
Then $f-v$ is smooth on $D$. Now $f\in{\mathcal F}^{p'}_{conf}(D,R)$ is evident for $p'\leq p$ and
$1<p'<\frac{4}{3}$.~\endproof


We now consider weak compactness property of a bounded sequence in ${\mathcal F}^p(D,R)$.
Let $\{f_k\}\subset\mathcal{F}^p(D,R)$. The blowup set of $\{f_k\}$ is defined to be
$$
\mathcal{C}(\{f_k\})=\left\{z\in D:\lim_{r\rightarrow 0}
\varliminf_{k\rightarrow+\infty}W(f_k,D_{r}(z))>\epsilon_2^2\right\}.
$$
Then for any $z\in D\backslash\mathcal{C}(\{f_k\})$,
we can find $r$ and a subsequence of $\{f_k\}$ which is still denoted by $\{f_k\}$ for simplicity, such that
$$\lim_{k\rightarrow+\infty}W(f_k,D_r(z))<\epsilon_0^2.$$
Then we get from Proposition \ref{epsilon} that
$\|f_k\|_{W^{2,p}(D_{r/2}(z))}<C\|\nabla f_k\|_{L^2(D_r)}$. Thus we may
assume $f_k$ converges weakly in $W^{2,p}(D\backslash
\mathcal{C}(\{f_k\}))$.

\begin{cor}\label{it}
Let $\{f_k\}\subset\mathcal{F}^p_{conf}(D,R)$ with
$$
\sup_k\{E(f_k,D)+W(f_k,D)\}<\Lambda<\infty
$$
and $f_0$ be the weak limit of $f_k$ in
$W^{2,p}_{loc}(D\backslash\mathcal{C}(\{f_k\}))$. If $p\in(1,\frac{4}{3})$,
then $f_0\in \mathcal{F}^p_{conf}(D,R)$ and
\begin{equation}\label{WL}
W(f_0,D)\leq\varliminf_{k\rightarrow+\infty}W(f_k,D).
\end{equation}
\end{cor}

\proof Set $\Delta f_k= F_k,k\in{\mathbb N}$. For any $\Omega\subset\subset D\backslash \mathcal{C}(f_k)$,
we have $\|f_k\|_{W^{2,p}(\Omega)}<C(\Omega)$. Then by the H\"older inequality and the Sobolev inequality
$$
\|F_k\|_{L^p(\Omega)}\leq\left\| \frac{1}{2}H(f_k)\nabla f_k\right\|_{L^2(\Omega)}\left\|\nabla f_k\right\|_{L^\frac{2p}{2-p}(\Omega)}
\leq C \Lambda^{\frac{1}{2}}\|\nabla f_k\|_{W^{1,p}(\Omega)}<C'(\Omega,\Lambda).
$$
We may assume, by selecting subsequences if necessary, that
$$
F_k\rightharpoonup F_0\s \hbox{in}\s L^{p}(\Omega)\s \hbox{and}\s
|H(f_k)||\nabla f_k|\rightharpoonup\alpha\s \hbox{in}\s L^{2}.
$$
Since  we may also assume
$|\nabla f_k|\rightarrow |\nabla f_0|$ in $L^2(\Omega)$ because $f_k\rightharpoonup f_0$ in $W^{2,p}(\Omega)$,
we have
$$|H(f_k)||\nabla f_k|^2\rightharpoonup \alpha|\nabla
f_0|$$
in the sense of measures in $\Omega$.
Define
$$\beta_0=\left\{\begin{array}{ll}
                    \frac{\alpha}{|\nabla f_0|}&
                     \hbox{ when}\, |\nabla f_0|\neq 0\\
                    0&\,\,\hbox{otherwise.}
                 \end{array}\right.$$
Clearly, $\beta_0|\nabla f_0|^2=\alpha|\nabla f_0|$. 
Let $F_k^+=\max\{F_k,0\}$ and $F_k^-=-\min\{F_k,0\}$. Then $F_k=F_k^+-F_k^-$ and $|F_k|=F_k^++F_k^-$.
We may assume that
$$
F_k^+\rightharpoonup F_0^1\s \hbox{and} \s F_k^-\rightharpoonup F_0^2\s \hbox{in}
\s L^p(\Omega).
$$
Obviously $F_0=F_0^1-F_0^2$. Then for any nonnegative function $\varphi\in C_0^\infty(\Omega)$,
$$
\int_\Omega\varphi |F_0|\leq \int_\Omega\varphi (F_0^1+F_0^2)
= \lim_{k\rightarrow+\infty}\int_\Omega
\varphi |F_k|\leq\lim_{k\rightarrow+\infty}
\int_\Omega \frac{1}{2}\varphi|H(f_k)||\nabla f_k|^2
=\int_\Omega\frac{1}{2}\varphi \beta_0|\nabla f_0|^2.
$$
Hence we conclude
$$
|F_0|\leq \frac{1}{2}\beta_0|\nabla f_0|^2, \s\hbox{a.e.}\, z\in D.
$$
Then, we have
$$
\int_\Omega\beta_0^2|\nabla f_0|^2\leq\int_\Omega\alpha^2
\leq\varliminf_{k\rightarrow+\infty}\int_\Omega|H_k(f_k)|^2|\nabla f_k|^2.
$$
Moreover, as $f_k$ converge in $L^2(\Omega)$, it follows from
$\partial f_k \otimes \partial f_k =0$ a.e. in $D$ that $\partial f_0\otimes \partial f_0=0$ a.e. in $D$ as $\Omega$ is arbitrary. Since $\sup_k\{E(f_k)+W(f_k)\}<\infty$, there are at most finitely many points in ${\mathcal C}(f_k)$.
Then we conclude that $f_0\in\mathcal{F}^p_{conf}(D,R)$ if $p\in(1,\frac{4}{3})$ by removing the point singularity
across $\mathcal{C}(f_k)$ ensured by Proposition \ref{rem}.
Furthermore, we have $H(f_0)\leq \beta_0$ whenever
$|\nabla f_0|\neq 0$, hence we get \eqref{WL}.
\endproof

\subsection{A criterion for absence of bubbles along cylinders}

Let $f_k\in\mathcal{F}^p_{conf}(Q(T_k),R)$, with
$$
\sup_k\{E(f_k)+W(f_k)\}<\Lambda<\infty.
$$
Given a sequence $t_k\in (-T_k,T_k)$ with
\begin{equation}\label{t}
T_k-t_k\rightarrow+\infty\,\,\,\hbox{and}\,\,\, t_k-(-T_k)
\rightarrow +\infty,
\end{equation}
we say the limit $f_0$ of a subsequence of $f_k(\theta,t+t_k)$, as in Corollary \ref{it}, is nontrivial if $E(f_0)>0$.   When $f_0$ is nontrivial, it is a bubble of $\{f_k\}$.

\begin{pro}\label{neck} Let $f_k\in{\mathcal F}^p_{conf}(S^1\times(-T_k,T_k),R)$ with
$$
\sup_k\{E(f_k)+W(f_k)\}=\Lambda<\infty.
$$
 Let $\epsilon_2$ be the constant in Proposition \ref{key}. If
$$
\lim_{T\rightarrow+\infty}\varliminf_{k\rightarrow+\infty}\sup_{t\in [-T_k+T,T_k-T]}W(f_k,S^1\times[t,t+1])<\epsilon_2^{2},
$$
then we have the following
\begin{enumerate}
\item  $\{f_k\}$ has no bubble;
\item there is no energy loss, i.e.
\begin{equation}\label{energyidentity}
\lim_{T\rightarrow+\infty}
\lim_{k\rightarrow+\infty}\int_{S^1\times[-T_k+T,T_k-T]}|\nabla f_k|^2=0,
\end{equation}
\item  there is no neck, i.e.
\begin{equation}\label{noneck}
\lim_{t\rightarrow+\infty}\lim_{k\rightarrow+\infty}
f_k(\theta,-T_k+t)=
\lim_{t\rightarrow+\infty}\lim_{k\rightarrow+\infty}
f_k(\theta,T_k-t).
\end{equation}
\end{enumerate}
\end{pro}

\proof  We may assume that $f_k(\theta,-T_k+t)$ and
$f_k(\theta,T_k-t)$ converge to $f_0^+(\theta,t)$  and
$f_0^-(\theta,t)$  weakly in $W^{2,p}_{loc}(S^1\times[0,+\infty))$, respectively.
Then $f_0^+(\phi)$ and $f_0^-(\phi)\in \mathcal{F}^p(D\backslash\{0\},R)$ with
$$
E(f^{\pm}_0(\phi))+W(f^{\pm}_0(\phi))\leq \Lambda
$$
 where $\phi$ is the conformal diffeomorphism between
$D\backslash\{0\}$ with $S^1\times(0,+\infty)$. By removability of
point singularities asserted in Proposition \ref{rem}, they are in $\mathcal{F}^{p'}(D,R)$ for some $p'>1$.
It then follows from the compact embedding $W^{2,p'}\subset L^2$:
$$
\lim_{T\to\infty}\int_{S^1\times [T,T+1]}\left(|\nabla f_0^+|^2+|\nabla f_0^-|^2\right)=0.
$$
Define
$f_k^*(t)=\frac{1}{2\pi}\int_0^{2\pi}f_k( t,\theta)d\theta$.
It is easy to check that
$$
\lim_{t\rightarrow+\infty}\lim_{k\rightarrow+\infty}\left|\int_{\partial
Q(T_k-t)}(f_k-f_k^*)\frac{\partial f_k}{\partial t}\right|=0.
$$
In fact, this can be seen as follows:
$$
\sup_{S^1\times\{T_k-T\}}|f_k-f_k^*|\leq C\, \underset{S^1\times\{T_k-T\}}{\hbox{osc}} f_k
$$
which will converge to $C \hbox{osc}_{S^1\times\{-T\}}f_0^+$ as $k\to\infty$.
By removability of singularity,
$$
\lim_{T\to\infty}\hbox{osc}_{S^1\times\{-T\}}f_0^+= 0.
$$
By the Sobolev trace embedding,
$$
\int_{S^1\times\{T_k-T\}}|\nabla f_k|\leq C\|\nabla f_k\|_{W^{1,p}(S^1\times[T_k-T-1,T_k-T+1])}
$$
By $\epsilon$-regularity,
$$
\left\|\nabla f_k\right\|_{W^{1,p}(S^1\times[T_k-T-1,T_k-T+1])}\leq C\|\nabla f_k\|_{L^2(S^1\times[T_k-T-2,T_k-T+2])}
<C.
$$
Then  \eqref{energyidentity} follows from \eqref{en}.

Let $m_k$ be the integer in $[T_k-T,T_k-T+1)$.
For $0\leq i\leq m_k-2$, applying Proposition \ref{key} on $S^1\times [i-m_k,m_k]$ (by shifting the center circle to $S^1\times \{i\}$, and the same below), we have
$$
\int_{S^1\times[i-2,i+2]}|\nabla f_k|^2<CE(f_k,Q(T_k-T))e^{-\delta(m_k-i)},\s
\delta=1-C\epsilon_2.
$$
Then from (\ref{f^*})
$$
\underset{S^1\times[i-1,i+1]}{\hbox{osc}}f_k\leq
C\sqrt{E(f_k,Q(T_k-T))}e^{-\frac{\delta}{2} (m_k-i)}.
$$
When $-m_k+2\leq i\leq 0$,  applying Proposition \ref{key} on
$S^1\times [-m_k,m_k+i]$, we get
$$
\int_{S^1\times[i-2,i+2]}|\nabla f_k|^2<CE(f_k,Q(T_k-T))e^{-\frac{\delta}{2}(m_k-|i|)},
$$
then we obtain
$$\underset{S^1\times[i-1,i+1]}{\hbox{osc}}f_k
\leq C\sqrt{E(f_k,Q(T_k-T))}e^{-\frac{\delta}{2} (m_k-|i|)}.$$
Hence,
$$
\underset{Q(T_k-T)}{\hbox{osc}}f_k\leq
2C\sqrt{E(f_k,Q(T_k-T))}\sum_{i=1}^{m_k}e^{-\frac{\delta}{2} (m_k-i)}\leq C'\sqrt{E(f_k,Q(T_k-T))}.
$$
Then \eqref{noneck} can be deduced from \eqref{energyidentity}.
\endproof

\subsection{Bubble trees for a sequence of maps from the disk $D$}
Let $f_k\in \mathcal{F}^p_{conf}(D,R)$ with
$$
\sup_k\left\{E(f_k,D)+W(f_k,D)\right\}=\Lambda<\infty.
$$
We assume
$0$ is the only blowup point of $\{f_k\}$, i.e. the only point such that
$$\lim_{r\rightarrow 0}\varliminf_{k\rightarrow +\infty}W(f_k,D_r(z))
\geq\epsilon_2^2.$$
We assume that
$f_k$ converges to $f_\infty$ weakly in $W^{2,p}_{loc}(D\backslash\{0\})$. The construction of  the bubble tree at $0$ will be divided into the following steps:

{\bf Step 1.} Construct the first level of the bubble tree.

There exists a sequence of points $z_k\in D$ and a sequence of radii $r_k\to 0$ such that
\begin{equation}\label{top}
W(f_k,D_{r_k}(z_k))=\frac{\epsilon_2^{2}}{2}
\end{equation}
and $W(f_k,D_r(z))<{\epsilon_2^2/2}$ for any $r<r_k$ and
$D_{r}(z)\subset D$. It is easy to check that $z_k\rightarrow 0$ as $0$ is the only blowup point of $\{f_k\}$.

We set $f_k'(z)=f_k(z_k+r_kz)$. Since
$\mathcal{C}(\{f_k'\})=\emptyset$,
$f_k'(z)$ converge weakly in $W^{2,p}_{loc}(\C)$.
We denote the limit by $f^F$. Note that  it may be a trivial mapping.

Let $(r,\theta)$ be the polar coordinates centered
at $z_k$, and set $T_k=-\ln r_k$. Let $\phi_k:S^1\times [0,T_k]\rightarrow\R^2$ be the conformal mapping
given by $\phi_k(t,\theta)=(e^{-t},\theta).$
Then
$$
\phi_k^*(dx^1\otimes dx^1+dx^2\otimes dx^2)=\frac{1}{r^2}(dt^2+d\theta^2).
$$
Thus $f_k\circ\phi_k\in \mathcal{F}^p_{conf}(S^1\times[0,T_k],R)$.
We will also denote $f_k\circ\phi_k$ by $f_k$ for simplicity of notations.

\begin{lem}\label{first-step} There exists a subsequence of $\{f_k\}$ and  $0=d_k^0<d_k^1<\cdots<d_k^l=T_k$ with $l<{\Lambda/\epsilon_2^2}+1$, such that
\begin{equation}\label{tree1}
\lim_{k\rightarrow+\infty}d_k^j-d_k^{j-1}=\infty,
\end{equation}
\begin{equation}\label{tree2}
W(f_k,S^1\times[d_k^j,d_k^j+1])\geq\epsilon_2^2, \s j\neq 0,l
\end{equation}
and
\begin{equation}\label{tree3}
\lim_{T\rightarrow+\infty}\varliminf_{k\rightarrow+\infty}\sup_{t\in[d_k^{j-1}+T,d_k^{j}-T]}W(f_k,S^1\times[t,t+1])\leq\epsilon_2^2,\s
 j=1,..., l.
\end{equation}
\end{lem}

\proof Suppose
$$(m-1)\epsilon_2^2<W(f_k,S^1\times[0,T_k])\leq\epsilon_2^2 \,m,$$
where $m$ is a positive integer. We prove the lemma by induction on $m$.

When $m=1$, the lemma is obvious by taking $d^0_k=0,d^1_k=T_k$ and (\ref{tree2}) is vacuous. Assuming the lemma
is true for $m-1$, we will prove it also true for $m$.
First of all, if \begin{equation}\label{ind}
\lim_{T\rightarrow+\infty}
\varliminf_{k\rightarrow+\infty}\sup_{t\in[T,T_k-T]}
W(f_k,S^1\times[t,t+1])\leq \epsilon_2^2,
\end{equation}
then the lemma follows since $[d^{j-1}_k+T,d^j_k-T]\subset[T,T_k-T]$. If \eqref{ind} does not hold,
we can find $t_k$ such that
$$t_k\rightarrow+\infty,\s T_k-t_k\rightarrow+\infty,$$
and
$$W(f_k,S^1\times[t_k,t_k+1])\geq\epsilon_2^2.$$
Then
$$
W(f_k,S^1\times[0,t_k])\leq\epsilon_2^2\,(m-1)\s \hbox{and}\s W(f_k,S^1\times[t_k+1,T_k])\leq\epsilon_2^2\,(m-1).
$$
Using the induction hypothesis on $[0,t_k]$ and $[t_k+1,T_k]$, we can find
$$
0=\bar{d}_k^0<\bar{d}_k^1<\cdots<\bar{d}_k^{\bar{l}}=t_k,\s
\hbox{and}\s t_k+1=\widehat{d}_k^0<\widehat{d}_k^1<\cdots<\widehat{d}_k^{\hat{l}}=T_k,
$$
such that
$$
\bar{d}_k^i-\bar{d}_k^{i-1}\rightarrow+\infty,\s  \s \widehat{d}_k^i-\widehat{d}_k^{i-1}\rightarrow+\infty,
$$
$$W(f_k,S^1\times[\bar{d}_k^j,\bar{d}_k^j+1])\geq\epsilon_2^2,\s
W(f_k,S^1\times[\widehat{d}_k^j,\widehat{d}_k^j+1])\geq\epsilon_2^2,$$
and
$$\lim_{T\rightarrow+\infty}\varliminf_{k\rightarrow+\infty}
\sup_{t\in[\bar{d}_k^{j-1}+T,\bar{d}_k^{j}-T]}
W(f_k,S^1\times[t,t+1])\leq\epsilon_2^2,$$
$$\lim_{T\rightarrow+\infty}\varliminf_{k\rightarrow+\infty}
\sup_{t\in[\widehat{d}_k^{j-1}+T,\widehat{d}_k^{j}-T]}
W(f_k,S^1\times[t,t+1])\leq\epsilon_2^2.$$
Put
$$
d_k^i=\left\{\begin{array}{ll}
                    \bar{d}_k^i&i\leq \bar{l},\\
                    \widehat{d}_k^{i-\bar{l}}&i>\bar{l}.
               \end{array}\right.
               $$
The induction is complete. \endproof

We now start to construct the bubble tree at the first level.
In Lemma \ref{first-step}, if $l=1$, in view of Proposition \ref{neck}, we do not do anything as there is no bubble
developing in $S^1\times [0,T_k]$ when $k\to\infty$.
If $l>1$, we set
$f_k^i(\theta,t)=f_k(\theta,d_k^i+t)$.
We may assume $\{f_k^i\}$ converges weakly in
$W^{2,p}$ to a bubble $f_\infty^i$ in any
compact set outside the blowup points of
$\{f_k^i\}$. By Proposition \ref{neck}, there is no other bubble of
$f_k$ between $f_\infty^i$
and $f_\infty^{i+1}$ and $f_\infty^i\cup f_\infty^{i+1}$
is connected.

Clearly, $\{f_k^0\}$ and $\{f_k^l\}$ have no blowup points. Moreover $f_\infty^0$ is just a part of $f_\infty$ and
$f_\infty^l$ is just a part of $f^F$. Removing the point singularity by Proposition \ref{rem},
$f_\infty^1$, $\cdots$, $f_\infty^{l-1}$ and $f^F$
can be considered as conformal mappings from
$S^2$ into $\R^n$.

\begin{center}
\includegraphics{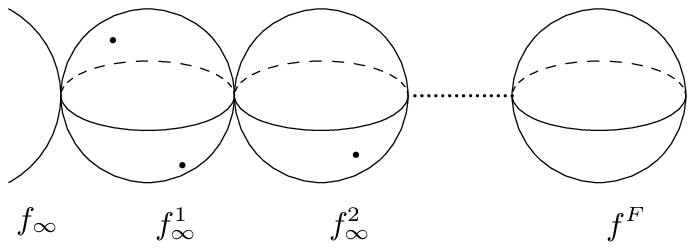}

\centering{Figure 2. Bubble tree: First level (dots denote concentration points)}
\end{center}

For a stratified sphere, we can define a {\it dual graph} as
following:
1) Associate one vertex for each component of the stratified sphere;
2) Vertices are connected by edges if the corresponding components meet at a point.

Let $S_1$ be the  stratified sphere with $l$ components whose dual graph is a tree, i.e. no loops.
We define $F^1$ to the continuous map from $S_1$ into $\R^n$, such that $F^1$ is  $f_\infty^i$ on the $i$-th
component  when $i<l$ and $f^F$  on the $l$-th component. We call $F^1$ {\it the first level of bubble
tree of $\{f_k\}$.}

We define $E(F^1)$ and $W(F^1)$ by
$$
E(F^1)=\sum_{i=1}^{l-1}\int_{S^1\times\R}|\nabla f_\infty^i|^2+\int_{S^2}|\nabla f^F|,\s
W(F^1)=\sum_{i=1}^{l-1}W(f_\infty^i)+W(f^F).
$$
Then
$$
\lim_{\delta\rightarrow 0}\lim_{k\rightarrow+\infty}\int_{D_\delta}
|\nabla f_k|^2=E(F^1)+\sum_i\sum_{p\in\mathcal{C}(\{f_k^i\})}\lim_{r\rightarrow 0}
\lim_{k\rightarrow+\infty}\int_{B_r(p)}|\nabla f_k^i|^2
$$
and
$$
\lim_{\delta\rightarrow 0}\lim_{k\rightarrow+\infty}W(f_k,D_\delta)\geq
W(F^1)+\sum_i\sum_{p\in\mathcal{C}(\{f_k^i\})}\lim_{r\rightarrow 0}
\lim_{k\rightarrow+\infty}W(f_k^i,B_r(p)).
$$

{\bf Step 2.}  We consider the convergence of $\{f_k^i\}$ near its blow up points.

For each $p\in \mathcal{C}(\{f_k^i\})$, we find a small $r$ such that $B_r(p)\subset S^1\times\R$ contains only one blowup point. Then for each $p$, using the arguments in Step 1, we will get the first level  of bubble tree of $\{f_k^i\}$,
which is a map $F_p$ from  a stratified sphere $S_p$  into $\R^n$.
Each $S_p$ is attached to $S_1$ at $p$. Taking union over $p\in \mathcal{C}(\{f_k^i\})$ gives us
a continuous map $F^2$ from $S_2$, which is a union of
stratified spheres, into $\R^n$. We call $F^2$  {\it  the second
level of the bubble tree of $\{f_k\}$.}

\begin{center}
\includegraphics{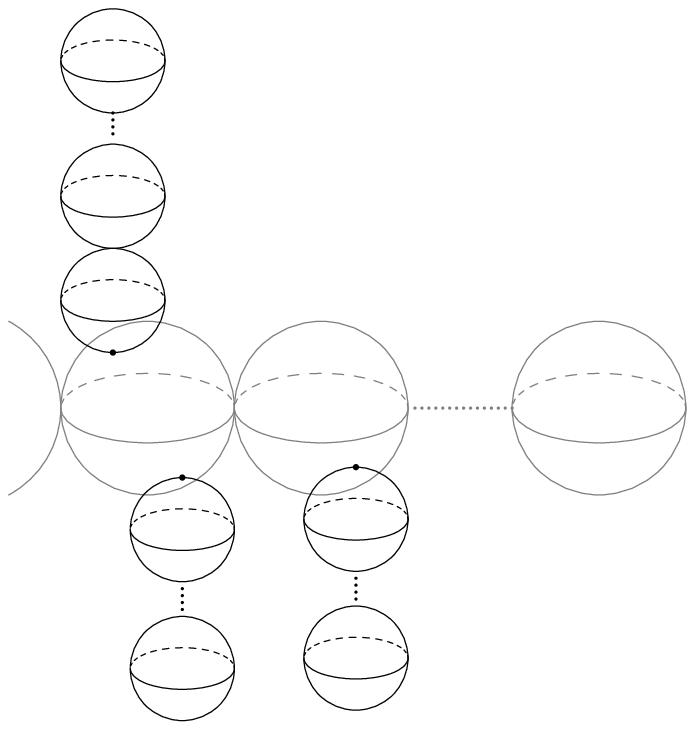}

\centering{Figure 3. Bubble tree: Second level}
\end{center}

{\bf Step 3.} In the same way, we can build the third and higher levels of the bubble tree.

Since each step will take away at least $\epsilon_2^2$ from the Willmore functional, the construction will stop after finite
many steps. In the end we get a stratified surface $S$ which is the union of all levels and a mapping $F$ from $S$ into $\R^n$. We shrink all the components of $S$
on which $F$ is trivial into points, {\it i.e.} deleting the ghost bubbles, then we get
a new stratified surface $S'$ and a continuous map $F'$
from $S'$ into $\R^n$, such that $F'$  is nontrivial on each component of $S'$. We call $F'$ is  {\it the bubble tree of
$\{f_k\}$ at $0$}. Moreover, we have
$$
\lim_{\delta\rightarrow 0}
\lim_{k\rightarrow+\infty}\int_{D_\delta}
|\nabla f_k|^2=E(F'),
$$
and
$$
W(F')\leq\lim_{\delta\rightarrow 0}\lim_{k\rightarrow+\infty}W(f_k,D_\delta).
$$

\subsection{Bubble trees for a sequence of maps from cylinders $Q(T_k)$
with $T_k\rightarrow+\infty$}
In this subsection, we show that in the situation that blowup occurs along a long cylinder we can divide the long cylinder into smaller ones and apply the results in the subsection 2.3 on each.

Let $f_k\in \mathcal{F}^p_{conf}(Q(T_k),R)$ with $T_k\rightarrow+\infty$. In light of Proposition \ref{neck}, we only need to consider the case that the following happens:
\begin{equation}\label{otherwise}
\lim_{T\rightarrow+\infty}\varliminf_{k\rightarrow+\infty}
\sup_{t\in [-T_k+T,T_k-T]}W(f_k,S^1\times[t,t+1])\geq \epsilon_2^{2}
\end{equation}
since otherwise there will be no bubbles, no necks and no energy loss. When (\ref{otherwise}) holds, there exist $t_k>0$ such that $T_k-t_k\to +\infty$ as $k\to\infty$ and
$$
W(f_k,S^1\times[t_k,t_k+1])\geq \epsilon^2_2.
$$

Case I. If  $\{t_k\}$ contains a bounded subsequence, we choose this subsequence and
 apply the bubble tree construction in subsection 2.3 at the blowup points.

Case II. If $\{t_k\}$ does not contain any bounded subsequences, then by Lemma \ref{first-step}, we can find (by translations)
$$
-T_k=d_k^0<d_k^1<\cdots<d_k^l=T_k
$$
which satisfy \eqref{tree1}, \eqref{tree2} and \eqref{tree3}. Recall that $l$ is independent of $k$.
We may assume $f_k^i(t,\theta)=f_k(d_k^i+t,\theta)$ converge
weakly to $f_\infty^i$ in $W^{2,p}$ outside the blowup points ${\mathcal C}(\{f^i_k\})$ of
$\{f_k^i\}$.
Let $\Sigma_\infty^1$ be the stratified surface with $l-1$ components whose dual
graph is a tree. Then we get a continuous map $F^1$ from $\Sigma^1_\infty$ into $\R^n$,
and $F^1$ is $f_\infty^i$ on the $i$-th component. Moreover, we have
$$
\lim_{T\rightarrow +\infty}\lim_{k\rightarrow+\infty}\int\limits_{S^1\times[-T_k+T,T_k-T]}
|\nabla f_k|^2=E(F^1)+\sum_{i=1}^{l-1}\sum_{p\in\mathcal{C}(\{f_k^i\})}
\lim_{r\rightarrow 0}\lim_{k\rightarrow+\infty}\int_{B_r(p)}|\nabla f_k^i|^2
$$
and
$$
\lim_{T\rightarrow +\infty}\lim_{k\rightarrow+\infty}W(f_k,Q(T_k-T))\geq
W(F^1)+\sum_{i=1}^{l-1}\sum_{p\in\mathcal{C}(\{f_k^i\})}\lim_{r\rightarrow 0}
\lim_{k\rightarrow+\infty}W(f_k^i,B_r(p)).
$$
The first level of the bubble tree of $\{f_k\}$ is  $F^1$ together with the bubble tree in Case I.
Then we repeat this process to construct the second level of the bubble tree at $\bigcup^l_{i=1}{\mathcal C}(\{f^i_k\})$, and similarly the third level and so on. The construction stops in finite steps.

\subsection{Convergence in Hausdorff distance}

The main aim of this subsection is to prove the following:
\begin{thm}\label{main1}
Assume that $(\Sigma,h_k)$ are a sequence of close  Riemann  surface of genus $g$,
where $h_k$ satisfies \eqref{metric-c}.
Suppose that $f_k\in \mathcal{F}^p_{conf}(\Sigma,h_k,R)$ with  $p\in(1,\frac{4}{3})$ and
\begin{equation}\label{energy.Willmore}
\sup_p\{E(f_k)+W(f_k)\}<\Lambda<\infty.
\end{equation}
Then either $f_k$ converges to a point, or
there is a stratified surface $\Sigma_\infty$ with $g(\Sigma_\infty)\leq g$, an $f_0\in \mathcal{F}^p_{conf}(\Sigma_\infty,R)$,
such that a subsequence of
$f_k(\Sigma_k)$ converges to  $f_0(\Sigma_\infty)$ in
the Hausdorff distance with
$$E(f_0)=\lim_{k\rightarrow+\infty}E(f_k),\s and\s W(f_0)\leq\lim_{k\rightarrow+\infty} W(f_k).$$
\end{thm}

\noindent{\bf Remark.} Here $f_0\in\mathcal{F}^p_{conf}(\Sigma_\infty,R)$ means that $f_0\in C^0(\Sigma_\infty,\R^n)$, and for any component $\Sigma^i_\infty$ of $\Sigma_\infty$, $F$ is
nontrivial on $\Sigma^i_\infty$ and $F|_{\Sigma^i_\infty}\in \mathcal{F}^p(\overline{\Sigma^i_\infty},h_i,R)$.\\

\noindent{\it Proof of Theorem \ref{main1}:} The proof will
be divided into three cases according to the genus of $\Sigma$.\vspace{0.7ex}

{\bf\noindent Spherical case.} When $\Sigma$ is a sphere, as there is only one conformal structure on a 2-sphere, we may let $h_k\equiv h$. Let $\mathcal{C}(\{f_k\})=\{p_1,\dots,p_m\}.$
We can choose $\delta$, such that $B_\delta(p_i)\cap B_\delta(p_j)=\emptyset$. Using isothermal coordinates,
 each $B_\delta(p_i)$ with metric $h$ is conformal to a Euclidean disk, the results can be deduced from subsection 2.3 directly. \vspace{0.7ex}

{\bf\noindent Toric case.} Suppose that $(\Sigma,h)$ is induced by lattice
$\{1,a+bi\}$ in $\C$, where $-\frac{1}{2}<a\leq\frac{1}{2}$, $b>0$, $a^2+b^2\geq 1$, and
$a\geq 0$ whenever $a^2+b^2=1$.
Then the conformal map $f$ from $(\Sigma,h)$ into $\R^n$ can be
lifted to a conformal map $\widetilde{f}$ from $\C$ into $\R^n$ which satisfies
$$\widetilde{f}(z+a+bi)=\widetilde{f}(z).$$
Let
$$\begin{array}{lll}
   \Pi: &\C\rightarrow S^1\times\R\\
        &a+bi\rightarrow (\pi a,\pi b),
\end{array}$$
be the conformal covering map. Then
$(\Sigma,h)$ is conformal to $(S^1\times\R)/G$,
where $G\cong \mathbb{Z}$ is the transformation group of
$S^1\times\R$ generalized by
$$
(\theta,t)\rightarrow (\theta+2\pi a,t+2\pi b).
$$
Then $f$ can be lifted to a map $f':S^1\times\R^n$, which satisfies $f'(\Pi)=\widetilde{f}$.

Now we assume $(\Sigma_k,h_k)=S^1\times\R/G_k$, where
$G_k$ is generalized by
$$(\theta,t)\rightarrow (\theta+\theta_k,t+b_k),\s where\s
b_k\geq \sqrt{\pi^2-\theta_k^2},\s and\s \theta_k\in [-\frac{\pi}{2},\frac{\pi}{2}].$$
In the moduli space $\mathcal{M}_1$ of genus 1 surfaces, $(\Sigma,h_k)$ diverges if and only if  $b_k\rightarrow+\infty$.

For $f_k\in \mathcal{F}^p_{conf}(\Sigma,h_k,R)$ with \eqref{energy.Willmore},
we lift each $f_k$ to a mapping $f_k':S^1\times\R\rightarrow\R^n$ which satisfies
$$f_k'(\theta,t)=f_k'(\theta+\theta_k,t+a_k).$$
After translations, we may assume that
$f_k'(\theta,t+\frac{a_k}{2})$ and $f_k'(\theta,
t-\frac{a_k}{2})$
have no blowup points as $k\to\infty$. Then $f_k'$ satisfies
the conditions in subsection 2.4 for $T_k=a_k/2$. Since
$f_k'(\theta,-T_k+t)=f_k'(\theta+\theta_k,T_k+t)$, the weak limit of
$f_k'(\theta,-T_k+t)$ in $W^{2,p}_{loc}(S^1\times[0,+\infty))$
and the weak limit of $f_k'(\theta,T_k+t)$ in $W^{2,p}_{loc}(S^1\times(-\infty,0])$
are just the two parts of a conformal map from $S^1\times\R$ into $\R^n$.
So the Hausdorff limit of $f_k(\Sigma)$ is the image of a continuous map $F$ from a
stratified surface  $S$ of genus 1 into $\R^n$ with
$$
E(F)=\lim_{k\rightarrow+\infty}E(f_k),\s
W(F)\leq \lim_{k\rightarrow+\infty}W(f_k).
$$

{\bf\noindent Hyperbolic case.}
For the hyperbolic case, we first briefly review the compactness of moduli space.

Let $\Sigma_0$ be a stable surface in $\overline{\mathcal{M}}_g$ with nodal points $\mathcal{N}=\{
a_1,\dots, a_{m'}\}$. Geometrically, $\Sigma_0$ is obtained by pinching
$m'$ non null homotopy curves in a surface with genus $g>1$ to points
$a_1,\dots,a_{m'}$,
thus $\Sigma_0\backslash\mathcal{N}$ can be divided
to finite components $\Sigma_0^1,\dots,\Sigma_0^s$. For each $\Sigma_0^i$, we can
extend $\Sigma_0^i$ to a smooth closed Riemann surface $\overline{\Sigma_0^i}$
by adding a point at each puncture. Moreover, the
complex structure of $\Sigma_0^i$ can be extended
smoothly to a complex structure of $\overline{\Sigma_0^i}$.

We say $h_0$ determines a hyperbolic structure on $\Sigma_0$ if $h_0$ is a  smooth complete metric on
$\Sigma_0\backslash\mathcal{N}$ with finite volume and Gauss curvature $-1$.
We define a neighborhood around each nodal point $a_j$ in $\Sigma_0$ by
$$
\Sigma_0(a_j,\delta)=\left\{p\in\Sigma_0:  \s \hbox{injrad}_{\Sigma_0\backslash\mathcal{N}}^{h_0}(p)<\delta,\s \forall p\in\Sigma_0(a_j,\delta)\backslash\{a_j\}\right\}\bigcup\,\{a_j\}.
$$
Let $h_0^i$ be the metric on $\overline{\Sigma_0^i}$ which has Gauss curvature $\pm1$ or
curvature 0, and is conformal to $h_0$ on $\Sigma_0^i$.

Now, we let $\Sigma_k$ be a sequence of compact Riemann surfaces of
fixed genus $g$ with hyperbolic structures $h_k$, such that $\Sigma_k
\rightarrow \Sigma_0$ in the moduli space
$\overline{\mathcal{M}_g}$. By Proposition 5.1 in \cite{Hum},
 there exists a maximal collection $\Gamma_k = \{\gamma_k^1,\ldots,\gamma_k^{m'}\}$
of pairwise disjoint, simple closed geodesics in $\Sigma_k$
with $\ell^j_k = L(\gamma_k^j) \to 0$, such that after passing
to a subsequence the following holds:
\begin{itemize}
\item[{\rm (1)}] There are maps $\varphi_k \in C^0(\Sigma_k,\Sigma_0)$,
such that $\varphi_k: \Sigma_k \backslash \Gamma_k \to \Sigma_0 \backslash \mathcal{N}$
is diffeomorphic and $\varphi_k(\gamma_k^j) = a_j$ for $j = 1,\ldots,m'$.
\item[{\rm (2)}] For the inverse diffeomorphisms
$\psi_k:\Sigma_0 \backslash \mathcal{N} \to \Sigma_k \backslash \Gamma_k$,
we have $\psi_k^\ast (h_k) \to h_0$ in $C^\infty_{loc}(\Sigma_0 \backslash
\mathcal{N})$, where $h_0$ determine a hyperbolic structure over
$\Sigma_0\backslash\mathcal{N}$.
\item[{\rm (3)}] Let $c_k$ be the complex structure over $\Sigma_k$, and $c_0$
be the complex structure over $\Sigma_0\backslash\mathcal{N}$. Then
$$\psi_k^*(c_k)\rightarrow c_0\s in\s
C^\infty_{loc}(\Sigma_0\backslash\mathcal{N}).$$
\end{itemize}

Moreover, we have the following Collar Lemma \cite{Halpern, Keen, M, Ra}:
\begin{lem}For each $\gamma_k^j$ as above, there is a collar $U_k^j$ containing $\gamma_k^j$, which is isometric to the cylinder $Q_k^j=Q(\frac{\pi^2}{l_k^j})$
with metric
\begin{equation}\label{metric}
h_k^j=\left(\frac{1}{2\pi\sin(\frac{l_k^j}{2\pi}t+
\theta_k)}\right)^2(dt^2+d\theta^2),
\end{equation}
where $\theta_k=\arctan(\sinh(
\frac{l_k^j}{2}))+\frac{\pi}{2}$.
Moreover, for any $(\theta,t)\in S^1\times
(-\frac{\pi^2}{l_k^j},\frac{\pi^2}{l_k^j})$, we have
\begin{equation}\label{injrad}
\sinh(\it{injrad}_{\Sigma_k}(\theta,t))\sin(
\frac{l_k^jt}{2\pi}+\theta_k)
=\sinh\frac{l_k^j}{2}.
\end{equation}
Let $\phi_k^j$ be the isometry between $Q_k^j$ and $U_k^j$. Then
$\varphi_k\circ\phi_k^{j}(\theta,\frac{\pi^2}{l_k^j}+t)\bigcup
\varphi_k\circ\phi_k^{j}(\theta,-\frac{\pi^2}{l_k^j}+t)$ converges in $C^\infty_{loc}(S^1\times(-\infty,0)\cup S^1\times(0,\infty))$ to an isometry from $S^1\times(-\infty,0)\cup S^1\times(0,+\infty)$ to $\Sigma_0(a_j,1)\backslash \{a_j\}$.
\end{lem}
The Collar Lemma can be found in \cite{Keen}, \cite{Halpern} and \cite{Hum}.

We also need the following local existence and compactness of conformal diffeomorphisms.

\begin{thm}\label{D.K.}\cite{D-K}
Let $h_k,h_0$ be smooth Riemannian metrics on a surface $M$,
such that $h_k \to h_0$ in $C^{s,\alpha}(M)$, where $s \in \N$,
$\alpha \in (0,1)$. Then for each point $z \in M$ there exist
neighborhoods $U_k, U_0$ of $z$ and smooth conformal diffeomorphisms
$\vartheta_k:D \to U_k, \vartheta_0:D\rightarrow U$, such that $\vartheta_k \to \vartheta_0$
in $C^{s+1,\alpha}(\overline{D},M)$.
\end{thm}

\noindent{\it Proof of Theorem \ref{main1} (continued):}
For a sequence $f_k\in\mathcal{F}^p_{conf}(\Sigma,h_k,R)$ satisfying the energy bound \eqref{energy.Willmore},
let
$$
\widetilde{f}_k=f_k\circ\psi_k
$$
which is a mapping from $\Sigma_0\backslash\mathcal{N}$ to $\R^n$. It is easy to check that $\widetilde{f}_k\in \mathcal{F}^p_{conf}(\Sigma_0\backslash\mathcal{N},\psi_k^{\ast}(h_k),R)$.

First, we show $\widetilde{f}_k$ converge in $W^{2,p}_{loc}(\Sigma_0\backslash(\mathcal{N}\cup \mathcal{C}(\{f_k\})))$. Given a point $z\in\Sigma_0\backslash(\mathcal{N}\cup\mathcal{C}(\{\widetilde{f}_k\}))$, we choose
 $U_k, U,\vartheta_k,\vartheta$ as in Theorem \ref{D.K.} and $U_k$, $U\subset \Sigma_0\backslash(\mathcal{N}\cup
\mathcal{C}(\{\widetilde{f}_k\}))$. Let
$$
\widehat{f}_k=\widetilde{f}_k\circ\varphi_k
$$
and note that
$\widehat{f}_k\in \mathcal{F}^p_{conf}(D,R)$. We can assume
that $\widehat{f}_k$ converge to $\widehat{f}_\infty$
in $W^{2,p}_{loc}(D_{3/4})$ with
$\partial \widehat{f}_\infty\otimes\partial \widehat{f}_\infty=0$.
Let $V=\vartheta(D_{1/2})$. Since $\vartheta_k$ converge to $\vartheta$, $\vartheta_k^{-1}(V)\subset D_{3/4}$ for sufficiently large $k$,
$\widetilde{f}_k=\widehat{f}_k(\vartheta_k^{-1})$ converge to $\widetilde{f}_\infty=\widehat{f}_\infty(\vartheta_0^{-1})$ weakly in $W^{2,p}(V,h_0)$. Then $\widetilde{f}_\infty
\in \mathcal{F}^{p}_{conf}(V,h_0,R)$.
Moreover, for any nonnegative function $\varphi$ with support in $V$, from Fatou's lemma
\begin{equation}\label{cut-off:W}
\lim_{k\rightarrow+\infty}\int_{V}\varphi |H(\widetilde{f}_k)|^2|\nabla \widetilde{f}_k|^2
=\lim_{k\rightarrow+\infty}\int_{D}\varphi(\vartheta_k)|H(\widehat{f}_k)|^2|
\nabla\widehat{f}_k|^2\geq \int_D\varphi(\vartheta_0)|H(\widehat{f}_0)|^2|\nabla\widehat{f}_0|^2.
\end{equation}
We may thus assume $\widetilde{f}_k$
converge weakly to $\widetilde{f}_\infty$ in $W^{2,p}_{loc}(\Sigma_0\backslash(\mathcal{N}\cup\,\mathcal
{C}(\{f_k\})))$.
Then $\widetilde{f}_\infty|_{\Sigma_0^i}\in
W^{2,p}_{loc}(\Sigma_0^i,h_0^i)$. So for $p\in (1,{4/3})$, $\widetilde{f}_\infty|_{\Sigma_0^i}$
can be extended to a map in $\mathcal{F}^p_{conf}(\overline{\Sigma_0^i},h_0^i,R)$.
Further,

$$
\lim_{k\rightarrow+\infty}E(f_k)=E(\widetilde{f}_\infty)+\sum_{z\in \mathcal{C}(\{\widetilde{f}_k\})}\lim_{r\rightarrow 0}\lim_{k\rightarrow+\infty}E(\widetilde{f}_k,B_r(z,h_0))
+\sum_j\lim_{\delta\rightarrow 0}\lim_{k\rightarrow+\infty}E(\widetilde{f}_k,\Sigma_0(a_j,\delta))
$$
and from (\ref{cut-off:W})
$$
\lim_{k\rightarrow+\infty}W(\widetilde{f}_k)\geq W(\widetilde{f}_\infty)+\sum_{z\in \mathcal{C}(\{\widetilde{f}_k\})}\lim_{r\rightarrow 0}\lim_{k\rightarrow+\infty}W(\widetilde{f}_k,B_r(z,h_0))+
\sum_j\lim_{\delta\rightarrow 0}\lim_{k\rightarrow+\infty}W(\widetilde{f}_k,\Sigma_0(a_j,\delta)).
$$

Next, we construct bubble trees at a point
$z\in\mathcal{C}(\{f_k\})\backslash\mathcal{N}$. We have
 a bubble tree $F$ of $\widehat{f}_k$ at $z$. We define it to be a  bubble tree of $\widetilde{f}_k$ at $z$.
By the arguments in subsection 2.3, we have
$$\lim_{r\rightarrow 0}\lim_{k\rightarrow+\infty}E(\widetilde{f}_k,B_r(z,h_0))
=\lim_{r\rightarrow 0}\lim_{k\rightarrow+\infty}E(\widehat{f}_k,D_r)=E(F),$$
and
$$\lim_{r\rightarrow 0}\lim_{k\rightarrow+\infty}W(\widetilde{f}_k,B_r(z,h_0))
=\lim_{r\rightarrow 0}\lim_{k\rightarrow+\infty}W(\widehat{f}_k,D_r)\geq W(F).$$

Lastly, we consider the convergence of $f_k$ at
the collars. Set $$\check{f}_k^j=f_k\circ\phi_k^j$$ and $T_k^j={\pi^2/l_k^j}-T$.
We may choose $T$ to be sufficiently large such that the two sequences
$\check{f}_k^j(T_k^j-t,\theta)$ and $\check{f}_k^j
(-T_k^j+t,\theta)$ have no blowup points in $[0,T]$ and $[-T,0]$ respectively (otherwise, return to the previous case as in Case I, subsection 2.4). Then $\check{f}_k^j$
satisfies the conditions in subsection 2.4. We get
a bubble tree $F^j$.
So the convergence of $\check{f}_k^j$ is clear. Since
$$
\check{f}_k^j=f_k\circ\phi_k^j=f_k\circ\psi_k\circ(\varphi_k\circ\phi_k^j)=
\widetilde{f}_k(\varphi_k\circ\phi_k^j),
$$
we have
\begin{eqnarray*}
\check{f}_k^j(T_k^j-t,\theta)&=&\widetilde{f}_k(\varphi_k\circ\phi_k^j(T_k^j-t,\theta)),\\
\check{f}_k^j(t-T_k^j,\theta)&=&\widetilde{f}_k(\varphi_k\circ\phi_k^j(t-T_k^j,\theta)).
\end{eqnarray*}
By the convergence statement of the Collar Lemma, $\varphi_k\circ\phi^j_k(t+T-\pi^2/l^j_k,\theta)$ and $\varphi_k\circ\phi_k^{j}({\pi^2/l_k^j}-T-t,\theta)$ converge in
$C^\infty_{loc}(S^1\times(0,\infty))$ to an isometry from $S^1\times(-\infty,0)\cup S^1\times(0,+\infty)$ to $\Sigma_0(a_j,1)\setminus \{a_j\}$.
We conclude that the image of the limit of $\check{f}_k^j(T_k^j-t,\theta)$ and that of $\check{f}_k^j(-T_k^j+t,\theta)$
are both contained in the image of $\widetilde{f}_\infty$.
Moreover,
$$\lim_{\delta\rightarrow 0}\lim_{k\rightarrow+\infty}E(\widetilde{f}_k,\Sigma_0(a_j,\delta))
=\lim_{T\rightarrow+\infty}\lim_{k\rightarrow+\infty}
E(\check{f},Q(T_k-T))=E(F^j),$$
and
$$\lim_{\delta\rightarrow 0}\lim_{k\rightarrow+\infty}W(\widetilde{f}_k,\Sigma_0(a_j,\delta))
=\lim_{T\rightarrow+\infty}\lim_{k\rightarrow+\infty}
W(\check{f},Q(T_k-T))=E(F^j).$$
Thus, we complete the proof.
\endproof

\begin{rem}\label{moduli} When $\Sigma_0\in \mathcal{M}_p$, i.e.
$\mathcal{N}=\emptyset$, $\psi_k$ is just a smooth
diffeomorphism
sequence from $\Sigma$ to $\Sigma$. In this case,
$g(\Sigma_\infty)=g(\Sigma)$, and
$$\Sigma_\infty=\Sigma_\infty'\cup S_1\cup S_2\cdots \cup S_m,$$
where $\Sigma_\infty'$ is a smooth Riemann surface of genus $p$, and each $S_i$ is a sphere.
\end{rem}

\begin{center}
\includegraphics{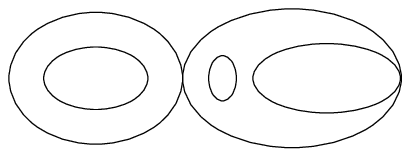}\hspace{5ex}\includegraphics{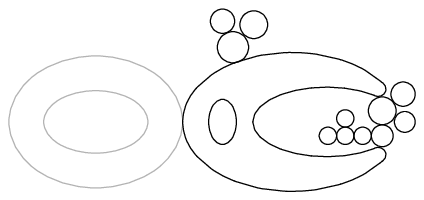}

\centering{Figure 4. $\Sigma_0$ (limit of $(\Sigma,h_k))$
and $\Sigma_\infty$} (some \\
components of $\Sigma_0$
do not appear in $\Sigma_\infty$)
\end{center}


We now generalize Theorem \ref{main1} to surfaces with marked points. Let us briefly review
the compactification of the moduli space of surfaces with marked points.
Let $\overline{\mathcal{M}}_{g,m}$ be the moduli space of closed
Riemann surfaces of genus $g$ with $m$ marked points.
Let $(\Sigma_0,x_{0,1},\dots,x_{0,m})\in\partial
\overline{\mathcal{M}}_{g,m}$ with nodal points $\mathcal{N}=\{
a_1,\dots, a_{m'}\}$. Geometrically, $\Sigma_0$ is obtained by pinching
some homotopically nontrivial closed curves which do not pass any of $x_{0,1},\dots,x_{0,m}$ into the points in ${\mathcal N}$, and $\Sigma\backslash\mathcal{N}$ can be divided
to connected components $\Sigma_0^1$, $\cdots$,
$\Sigma_0^s$. For each $\Sigma_0^i$, we can
extend $\Sigma_0^i$ to a smooth closed Riemann surface $\overline{\Sigma_0^i}$
by adding a point at each puncture. Moreover, the
complex structure of $\Sigma_0^i$ can be extended
smoothly to a complex structure of $\overline{\Sigma^i_0}$.

We say  $h$  is a hyperbolic structure on
$(\Sigma,x_1,\dots,x_m)\in\mathcal{M}_{g,m}$ if $h$ is a smooth complete metric on $\Sigma\backslash\{x_1,\dots,x_m\}$
with curvature $-1$ and finite volume.
We say  $h_0$  is a hyperbolic structure on
$(\Sigma_0,x_{0,1},\dots,x_{0,m})\in\overline{\mathcal{M}}_{g,m}\backslash
\mathcal{M}_{g,m}$ if $h_0$ is a smooth complete metric on $$\Sigma\backslash\{a_1,\dots,a_{m'},x_{0,1},...,x_{0,m}\}$$
with curvature $-1$ and finite volume.

For a surface $\Sigma$ with hyperbolic structure $h$ and with marked points $x_1,\dots,x_m$, we define
$\Sigma^*=\Sigma\backslash\{x_1,\dots,x_m\}$, and $h^*$ to be the hyperbolic structure on $(\Sigma,x_1,\dots,x_m)$ which is conformal to $h$
on $\Sigma^*$.

Let $\{(\Sigma_k,x_{k,1},\dots,x_{k,m}\}$ be a sequence of marked surfaces in $\mathcal{M}_{g,m}$
with hyperbolic structures $h_k$ and
$$
(\Sigma_k,x_{k,1},\dots,x_{k,m})\rightarrow (\Sigma_0,x_{0,1}, \dots,x_{0,m})\s in\s \overline{\mathcal{M}}_{g,m}.
$$
By Proposition 5.1 in \cite{Hum} again, there exists a maximal collection $\Gamma_k = \{\gamma_k^1,\ldots,\gamma_k^{m'}\}$ of pairwise disjoint, simple closed geodesics in $\Sigma_k$
with $\ell^j_k = L(\gamma_k^j) \to 0$ as $k\to\infty$, such that, after passing to a subsequence if necessary, the following holds:
\begin{itemize}
\item[{\rm (1)}] There are maps $\varphi_k \in C^0(\Sigma_k,\Sigma_0)$, such that $\varphi_k: \Sigma_k \backslash \Gamma_k \to \Sigma_0 \backslash \mathcal{N}$
is diffeomorphic and $\varphi_k(\gamma_k^i) = a_i$ for $i = 1,\dots,m'$, and $\varphi_k(x_{k,j}) = x_{0,j}$ for $j = 1,\dots,m$.
\item[{\rm (2)}] For the inverse diffeomorphisms
$\psi_k:\Sigma_0 \backslash \mathcal{N} \to \Sigma_k \backslash \Gamma_k$,
we have $\psi_k^\ast (h_k) \to h$ in $C^\infty_{loc}(\Sigma_0^* \backslash
\mathcal{N})$.
\item[{\rm (3)}] Let $c_k$ be the complex structure over $\Sigma_k$, and $c_0$
be the complex structure over $\Sigma_0\backslash\mathcal{N}$. Then
$$\psi_*(c_k)\rightarrow c_0\s in\s
C^\infty_{loc}(\Sigma_0\backslash\mathcal{N}).$$
\end{itemize}

Moreover, the Collar Lemma also holds for the moduli space with marked points.

\begin{thm}\label{main2}
In addition to the assumptions in Theorem \ref{main1}, we assume for $m\geq 2$
$$
 y_1, \dots ,y_m\in f_k(\Sigma).
$$
Then there is a stratified surface $\Sigma_\infty$ with $g(\Sigma_\infty)\leq g$, and an $f_0\in \mathcal{F}^p_{conf}(\Sigma_\infty,R)$ with
$$
y_1,\dots,y_m\in f_0(\Sigma_\infty),
$$
such that a subsequence of
$f_k(\Sigma_k)$ converges to  $f_0(\Sigma_\infty)$ in Hausdorff distance
with
$$E(f_0)=\lim_{k\rightarrow+\infty}E(f_k)\s and\s W(f_0)\leq\lim_{k\rightarrow+\infty} W(f_k).$$
\end{thm}

\proof


Let $\widetilde{f}_k=f_k\circ\psi_k$. In view of Theorem \ref{main1},
we only need to consider convergence of $\{\widetilde{f}_k\}$ near $x_{0,j}, j=1,\dots,m$.

Choose a complex coordinate  $\{U,(x,y)\}$ on $\Sigma_0$  compatible with $c_0$,
 with $x_{0,j}=(0,0)$. Let $c_k'=\psi_k^{\ast}(c_k)$.
We set
$$e_1=\frac{\partial }{\partial x},\s e_2=c_k'(e_1),$$
and $h_k'$ to be the metric on $U$ defined by
$$h_k'(e_1,e_1)=h_k'(e_2,e_2)=1,\s h_k'(e_1,e_2)=0.$$
Then $h_k'$ is compatible  with $c_k'$, and converges
smoothly to a metric which is compatible  with $c_0$
in $U$. Then we consider the weak convergence of $\{\widetilde{f}_k\}$
in $U\backslash\mathcal{C}(\{\widetilde{f}_k\})$, using the arguments
in subsection 2.3.

It remains to check that each marked point $y_i$ is on the image of $f_\infty$ or one of the bubbles. If $x_{0,j}$ is not a blow up point of $\{\widetilde{f}_k\}$, it is obvious that $y_j\in f_\infty(U)$. Now assume $x_{0,j}$ is the only blow-up point in $D$. We take $U_k$, $U_0$, $\vartheta_k,\vartheta_0,\widehat{f}_k,V$ as in the proof of Theorem \ref{main1} for $z=x_{0,j}$.
We will prove it by induction on the number of the levels of the bubble tree.
We take $z_k$, $r_k$, $\phi_k$ and $d_k$ as in subsection 2.3 for $\widehat{f}_k$.  If ${|z_k|/r_k}<L$ for some fixed $L$, then we may assume $-{z_k/r_k}\rightarrow z_\infty$, by selecting a subsequence if necessary. Recalling that $\widehat{f}_k(0)\equiv y_j$, we get $y_j=\widehat{f}^F(z_\infty)$. Let $(r,\theta)$ be the polar coordinates centered
at $z_k$,  $T_k=-\ln r_k$ and $\phi_k:[0,T_k]\times S^1\rightarrow\R^2$ be the conformal mapping
given by $\phi_k(t,\theta)=(e^{-t},\theta)$.
We set $\phi_k^{-1}(0)=(t_k,\theta_k)$.
Then
${|z_k|/r_k}\rightarrow+\infty$ means that $t_k\in [0,T_k]$ and
$T_k-t_k\rightarrow+\infty$.
Thus we may assume $t_k\in [d_k^i,d_k^{i+1}]$ for some $i$,
where $d_k^i$ are defined in Lemma \ref{first-step}. Then, if $t_k-d_k^i\rightarrow+\infty$
and $d_k^{i+1}-t_k\rightarrow+\infty$, we have
$y_j=f_\infty^i(+\infty)=f_\infty^{i+1}(-\infty).$
If at least one of $t_k-d_k^i$
and $d_k^{i+1}-t_k$ is bounded above for all large $k$, then we repeat the above argument at the second level of the bubble tree, and proceed in this way for the finitely many levels of the bubble tree if necessary, and we conclude  that $y^j$ is on one of bubbles of $\widetilde{f}_k^i$ or $\widetilde{f}_k^{i+1}$.

Finally, as $m\geq 2$ and all $y_i\in f_0(\Sigma_\infty)$, $f_k$ cannot converge to a single point.      \endproof

\section{Branched conformal immersions and proof of
Theorem \ref{main}}

For a branched conformal immersion, we have the following result:

\begin{thm}\label{removal}\cite{K-L}
Suppose that $f\in W^{2,2}_{conf,loc}(D\backslash \{0\},\R^n)$ satisfies
$$
\int_D |A_f|^2\,d\mu_g < \infty \quad \mbox{ and } \quad \mu_g(D) < \infty,
$$
where $g_{ij} = e^{2u} \delta_{ij}$ is the induced metric. Then
$f \in W^{2,2}(D,\R^n)$ and we have
\begin{eqnarray*}
u(z) & = & \lambda\log |z|+ \omega(z) \quad \mbox{ where }
\lambda\geq 0, \quad\lambda \in \mathbb{Z},  \quad\omega \in C^0 \cap W^{1,2}(D),\\
-\Delta u & = & -2\lambda\pi \delta_0+K_g e^{2u} \quad \mbox{ in }D.
\end{eqnarray*}
The density  of $f(D_\sigma)$ as  varifolds  at $f(0)$ is given by
$\lambda+1$ for any small
$\sigma > 0$.
\end{thm}

The classical Gauss-Bonnet formula is generalized in \cite{E-T}
for smooth branched surface. Following arguments in \cite{E-T}, we provide a
version for $W^{2,2}$ branched conformal immersions.

\begin{lem}\label{GB} Let $(\Sigma,g)$ be a smooth closed compact Riemann surface. Then
for any $f\in W_{b,c}^{2,2}(\Sigma,g,\R^n)$, there holds
\begin{equation}\label{G.B.}
\int_\Sigma K_fd\mu_f=2\pi\chi(\Sigma)+2\pi b,
\end{equation}
where $b$ is the number of branch points counted with
multiplicities and at each branch point $p$ the branching order is $\lambda=\theta^2(p)-1$.
\end{lem}

\proof Without loss of generality, we only prove the case that $f$
has only one branch point $p$. Let $g_f=e^{2u}g$ be the metric induced by
$f$ and $K_f$ be its Gauss curvature. In \cite{K-L}, we proved that the equation
$-\Delta_gu=K_fe^{2u}-K_g$ holds weakly in $\Sigma\backslash\{p\}$: for any smooth
$\varphi$ with support in $\Sigma\backslash\{p\}$, it holds
$$
\int_\Sigma\nabla_g u\nabla_g\varphi \,d\mu_g=
\int_\Sigma\varphi K_fe^{2u}d\mu_g-\int_\Sigma\varphi K_g du_g.
$$

Take a complex coordinate chart $\{U;z\}$
around $p$ with $p=0$. For any small $\epsilon>0$, we
choose a function $\varphi_\epsilon(z)=\varphi_\epsilon(|z|)$ between 0 and 1
with $|\varphi_\epsilon'|<{C/\epsilon}$ and equals 1 outside
$D_\epsilon$ and $0$ in $D_{\epsilon/2}$. Then we have
$$
\int_{D_\epsilon}\frac{\partial u}{\partial r}
\varphi_\epsilon' dx=
\int_\Sigma\varphi_\epsilon K_fe^{2u}d\mu_g-\int_\Sigma\varphi_\epsilon K_g
du_g.
$$
By Theorem \ref{removal}, $u=\lambda\log|z|+\omega$, we get
$$
\int_{D_\epsilon}\frac{\partial u}{\partial r}
\varphi_\epsilon' dx=\int_{D_\epsilon}\frac{\partial \omega}{\partial r}\varphi_\epsilon'dx+2\pi\lambda\left(\varphi_\epsilon(\epsilon)-\varphi_\epsilon(0)\right)=
\int_{D_\epsilon}\frac{\partial \omega}{\partial r}\varphi_\epsilon'dx+2\pi\lambda.
$$
Since
$$\int_{D_\epsilon}\left|\frac{\partial \omega}{\partial r}\varphi_\epsilon'\right|
\leq C\left(\int_{D_\epsilon\backslash D_{\epsilon/2}}\left|\frac{\partial\omega}{\partial r}\right|^2\right)^{1/2}
\left(\int_{D_\epsilon\backslash D_{\epsilon/2}}\frac{1}{r^2}\right)^{1/2}\leq C \|\nabla\omega\|_{L^2(D_\epsilon)}
\rightarrow 0,\s \hbox{as}\s \epsilon\rightarrow 0,
$$
we conclude, by applying the classical Gauss-Bonnet theorem on $(\Sigma,g)$, that
$$
\int_{\Sigma}K_fd\mu_f=\lim_{\epsilon\rightarrow 0}
\int_{\Sigma}\varphi_\epsilon K_fd\mu_f=
\lim_{\epsilon\rightarrow 0}
\int_{\Sigma}\varphi_\epsilon K_gd\mu_g+\lim_{\epsilon\rightarrow 0}
\int_{\Sigma}\frac{\partial u}{\partial r}\varphi_\epsilon'd\mu_g=2\pi\chi(\Sigma)+2\lambda\pi
$$
and complete the proof.
\endproof

\begin{rem}\label{rem-gb} Since $\int_\Sigma K_fd\mu_f\leq W(f)$, it follows from Lemma \ref{GB}
$$b\leq \frac{1}{2\pi}W(f)-\chi(\Sigma).$$
Moreover,
$$\int_\Sigma|A_f|^2d\mu_f=4W(f)-2\int K_f\leq 4 W(f)-2\pi\chi(\Sigma).$$
Then $\sup_kW(f_k)<+\infty$ implies that $\sup_kb_k<+\infty$ and $\sup_k\int_\Sigma|A_{f_k}|^2d\mu_{f_k}<+\infty$.
\end{rem}

To study the convergence conformal immersions,
we recall an important result of H\'elein.

\begin{thm}\label{Helein} \cite{He} Let $f_k\in W^{2,2}_{conf}(D,\R^n)$
be a sequence of conformal immersions with induced metrics
$(g_k)_{ij} = e^{2u_k} \delta_{ij}$, and assume
$$
\int_D |A_{f_k}|^2\,d\mu_{g_k} \leq \gamma <
\gamma_n =
\begin{cases}
8\pi & \mbox{ for } n = 3,\\
4\pi & \mbox{ for }n \geq 4.
\end{cases}
$$
Assume also that $\mu_{g_k}(D) \leq C$ and $f_k(0) = 0$.
Then $f_k$ is bounded in $W^{2,2}_{loc}(D,\R^n)$, and there
is a subsequence such that one of the following two alternatives
holds:
\begin{itemize}
\item[{\rm (a)}] $u_k$ is bounded in $L^\infty_{loc}(D)$ and
$f_k$ converges weakly in $W^{2,2}_{loc}(D,\R^n)$ to a conformal
immersion $f \in W^{2,2}_{conf,loc}(D,\R^n)$.
\item[{\rm (b)}] $u_k \to - \infty$ and $f_k \to 0$ locally uniformly on $D$.
\end{itemize}
\end{thm}

H\'elein first proved the above result  for $\gamma_n=
{8\pi/3}$ \cite[Theorem 5.1.1]{He}. In \cite{K-L} $\gamma_n$ is shown to be optimal.

Before proving Theorem \ref{main}, we recall a monotonicity formula (for more
details, see \cite{K-S2, S}).
Let $\mu$ be a 2-dimensional integral varifold
with square integrable weak mean curvature
$H_\mu\in L^2(\mu)$. Then we have
$$
g_{x_0}(\varrho) \leq g_{x_0}(\sigma), \s \hbox{when}\s  \varrho<\sigma,
$$
where
$$
g_{x_0}(r) =
\frac{\mu(B_r(x_0))}{\pi r^2}
+ \frac{1}{4 \pi} W(\mu,B_r(x_0))
+ \frac{1}{2\pi r^2} \int_{B_r(x_0)} \langle x-x_0,H \rangle \,d\mu.
$$

When $\Sigma$ is compact
and connected, if we let $\sigma\rightarrow+\infty$, and $\varrho\rightarrow 0$, then the area density of $\Sigma$ at $x_0$ satisfies
\begin{equation}\label{density1}
\theta^2(\mu,x_0) \leq \frac{1}{4\pi} W(f).
\end{equation}
If we only let $\varrho\rightarrow 0$, then we get
\begin{equation}\label{density2}
\theta^2(\mu,x_0) \leq \frac{\mu(B_\sigma(x_0))}{\pi \sigma^2}+ CW(f,B_\sigma(x_0)) +C\left( \frac{\mu(B_\sigma(x_0))}{\pi \sigma^2}  \right)^{\frac{1}{2}}W(f,B_\sigma(x_0))^{\frac{1}{2}}.
\end{equation}
Another useful consequence (cf. \cite{S}, \cite{K-S2}) is the following:
\begin{equation}\label{diameter}
\Big(\frac{\mu(\Sigma)}{W(f)}\Big)^{\frac{1}{2}}
\leq \diam f(\Sigma) \leq
C \Big(\mu(\Sigma)\, W(f)\Big)^{\frac{1}{2}}.
\end{equation}

\noindent{\it Proof of  Theorem \ref{main}.}
Consider a branched conformal immersion $f_k
\in W^{2,2}_{b,c}(\Sigma,h_k,\R^n)$, where $h_k$
satisfies \eqref{metric-c}. The following equation clearly holds on $\Sigma$ away from the finitely many branch points; the singularities at the branch points can be removed by Theorem \ref{removal}, thus it holds on entire $\Sigma$:
$$\Delta_{h_k} f_k=\frac{1}{2}H_{f_k}|\nabla_{h_k}f_k|^2.$$
By Remark \ref{rem-gb}, the number of branch points and
$\|A_{f_k}\|_{L^2}$ are both bounded from above.

By \eqref{diameter},
 $\mbox{\rm diam} f_k(\Sigma)\leq{R}$ for some $R>0$.
 Then $f_k\in \mathcal{F}^2_{conf}(\Sigma,h_k,R+R_0)$. We only need
 to prove that $f_0$ in Theorem \ref{main1} and Theorem \ref{main2}
is also branched conformal  immersion.

In fact, we only need to prove the following:
If $f_k$ are branched conformal immersions from
$D$ into $\R^n$ with a uniform upper bound on the number of branch points, then  the limit
$f_0$ is either a point or a branched conformal immersion.
Let $P$ be the limit set of the branch points,
and
$$\mathcal{S}(\{f_k\})=\{z\in D:\lim_{r\rightarrow 0}
\varliminf_{k\rightarrow+\infty}\int_{D_r(z)}|A_k|^2\geq
\hat{\epsilon}^2\},$$
where $\hat{\epsilon}\leq\min\{\sqrt{4\pi},4\epsilon_0\}$.
By Theorem \ref{Helein}, $f_k$ will converge weakly in $W^{2,2}_{loc}
(D\backslash(\mathcal{S}\cup P))$ to either a conformal
immersion or a point. If the limit is not a single point, by Theorem \ref{removal},
the limit can be extended across the finite set ${\mathcal S}\cup P$ to a branched conformal immersion of $D$, hence $\Sigma$, in $\R^n$.
\endproof

\section{Willmore functional for  surfaces in compact manifolds}

Let $N$ be a compact Riemannian manifold without boundary.
We embed $N$ into $\R^n$ isometrically so that any immersion  of $\Sigma$ in $N$
can be regarded as an immersion in $\R^n$. Let $A_{\Sigma,N},A_{\Sigma,\R^n}$
and $A_{N,\R^n}$ be the second fundamental forms of $\Sigma$ in $N$, in $\R^n$ and $N$ in $\R^n$ respectively. The $L^2$ integrals of these quantities can be related as in the following simple lemma.

\begin{lem}\label{H}
For any $f\in W^{2,2}_{b,c}(\Sigma,h,N)$, we have
\begin{equation}\label{4.11}
\int_{\Sigma}|H_{f,\Sigma,\R^n}|^2d\mu_{f}
\leq C\mu(f)+\int_\Sigma|H_{f,\Sigma,N}|^2d\mu_f,
\end{equation}
and
\begin{equation}\label{4.12}
\int_{\Sigma}|A_{f,\Sigma,\R^n}|^2d\mu_{f}\leq C\int_{\Sigma}(1+|H_{f,\Sigma,N}|^2)d\mu_{f}+C',
\end{equation}
where $C$ only depends on $N$ and $C'$ only depends on the Euler characteristic of $\Sigma$.
\end{lem}

\proof Let $e_1,\dots,e_n$ be an orthonormal frame of $T_{f_k(p)}\R^n$ with
$e_1$, $e_2\in Tf(\Sigma)$ and $e_1,\dots,e_k\in TN$. Then for $f_i=\frac{\partial f}{\partial x^i}$ we have

$$
\nabla^N_{f_i} f_j=\sum_{l=1}^k\lan f_{ij},e_l\ran e_l
$$
and
$$
A_{\Sigma,N}(f_i,f_j)=\sum_{m=3}^k\lan\nabla^N_{f_i}f_j,
e_m\ran e_m=\sum_{m=3}^k\lan f_{ij},e_m\ran e_m.
$$
Thus if $F = i \circ f $ where $i:N\to\R^n$ is the isometric embedding, we have
$$
A_{\Sigma,\R^n}(F_i,F_j)=\sum_{m=3}^n\lan F_{ij},e_m\ran e_m
=A_{\Sigma,N}(f_i,f_j)+A_{N,\R^n}(F_i,F_j).$$
Hence,  we have
\begin{equation*}
H_{\Sigma,\R^n}(f)=H_{\Sigma,N}(f)+g^{ij}A_{N,\R^n}(F_i,F_j).
\end{equation*}
Noting that $H_{\Sigma,N}(f)\perp A_{N,\R^n}$, we get
$$
\left|H_{\Sigma,\R^n}(f)\right|^2=\left|H_{\Sigma,N}(f)\right|^2+
\left|g^{ij}A_{N,\R^n}(F_i,F_j)\right|^2\leq \left|H_{\Sigma,\R^n}\right|^2+\|A_{N,\R^n}\|_{L^\infty}
$$
where $\|A_{N,\R^n}\|_{L^\infty}$ is bounded since $N$ is compact.
Thus, integrating over $\Sigma$ yields  \eqref{4.11}, and by \eqref{4.11} and Remark \ref{rem-gb}, we get
$$
\int_\Sigma|A_{\Sigma,\R^n}(f)|^2d\mu_{f}<C\left(1+\mu(f)+\int_\Sigma|H_{\Sigma,N}(f)|^2d\mu_{f}\right).
$$

\endproof

\subsection{Willmore sphere passing through fixed points}
In this subsection, we let
$$
W_n(f)=\int_{S^2}\left(1+\frac{1}{4}\left|H_f\right|^2\right)d\mu_f
$$
where $f$ is a $W^{2,2}$ conformal immersion of $S^2$ in the round unit sphere ${\mathbb S}^n$ for some $n>2$.
We consider the existence of minimizers of
$$
\beta_0^n(y_1,\dots,y_m)=\inf \{W_n(f):{y_1,\dots,y_m\in f(S^2)}\}
$$
where $y_1,\dots,y_m$ are fixed distinct points in ${\mathbb S}^n$. When $m\geq 2$, $\beta^n_0(y_1,\dots,y_m)$ is positive by the conformality of the functions $f$.

\begin{pro}\label{existence:sphere} When $m\geq 2$, any $\beta_0^n(y_1,\dots,y_m)<8\pi$ is attained
by a $W^{2,2}$-conformally embedded $S^2$ in ${\mathbb S}^n$.
\end{pro}

\proof Let $\{f_k\}$ be a minimizing sequence
of $\beta_0^n(y_1,\dots,y_m)$. We can consider
$f_k$ as conformal map from $S^2$ into $\R^n$. By
Theorem \ref{main}, $f_k$ will converge to a mapping $f_0$ which is a $W^{2,2}$
branched conformal immersion from a stratified
sphere $\Sigma_\infty$ into ${\mathbb S}^n$ with
$$
y_1,\dots,y_m\in f_0(\Sigma_\infty),\s W_n(f_0)\leq\beta_0^n(y_1,\dots,y_m)<8\pi.
$$
Composing with a stereographic projection $\Pi$ from ${\mathbb S}^n$ minus a point not on $f_0(S^2)$ into $\R^n$,
we see $W_n(f_0)=W(\Pi\circ f_0)$ and $\theta^2_{f_0(p)}=\theta^2_{\Pi\circ f(p)}$. Now, by \eqref{density1} we have
$$
\theta^2_{f_0(p)}\leq\frac{1}{4\pi} W_n(f_0).
$$
By Theorem \ref{removal}
$$
\lambda(p)+1=\theta^2_{f(p)}\leq \frac{1}{4\pi}W_n(f_0)<2
$$
thus $\lambda(p)=0$ which means $f_0$ has no branched points. Moreover, that the area density of $\Sigma_\infty$ is one everywhere implies that
$\Sigma_\infty$ has only 1 component and
$f_0$ has no intersection points.
Thus $\Sigma_\infty=S^2$, and $f_0$ is an (Lipschitz) embedding.
\endproof

\begin{cor}
For any $\epsilon>0$, there is a Willmore sphere $f:S^2\to{\mathbb S}^n$
 with $W_n(f)<4\pi+\epsilon$, which has
at least 2 nonremovable singular points.
\end{cor}

\proof
Take five distinct points $y_1,\dots,y_5\in {\mathbb S}^n$,
such that there is no round 2-sphere passing through all of them. Recall the Willmore functional $W_n$ of a round 2-sphere is $4\pi$. We can choose the five points to be very closed to a
round 2-sphere, such that there is a 2-sphere $\Sigma$ which is not round and contains
$y_1,\dots,y_5$ with
$$
W_n(\Sigma)<4\pi+\epsilon.
$$
Then we can find a $W^{2,2}$ conformal embedding  $f:S^2\rightarrow {\mathbb S}^n$, such that  $f(S^2)$ passes  through $y_1,\dots,y_5$, and attains $\beta_0^n(y_1,\dots,y_5)$, by Proposition \ref{existence:sphere}.

Choose a point $P\in {\mathbb S}^n\backslash\Sigma$ as the north pole.
Let $\Pi$ be the stereographic projection from
${\mathbb S}^n\backslash\{P\}$ to $\R^n$, and denote $\widetilde{y}_i=\Pi(y_i)$ and
$\widetilde{f}=\Pi(f)$. By the conformal invariance of the
Willmore functional, we have
$$
W_n(f)=\frac{1}{4}\int_{S^2}|H_{\widetilde{f}}|^2d\mu_{\widetilde{f}}.
$$
Then $\widetilde{f}$ attains
$$\inf\left\{\frac{1}{4}\int_{S^2}|H_{\varphi}|^2d\mu_\varphi:
\varphi\in W^{2,2}_{conf}(S^2,\R^n),\s \widetilde{y}_1,\dots,\widetilde{y}_5\in \varphi(S^2)\right\}.$$
Then by  results in \cite{R}, $\widetilde{f}(S^2)$ is smooth on $\widetilde{f}(S^2)\backslash\{\widetilde{y}_1,\dots,\widetilde{y}_5\}$.
However, the Gap Lemma
in \cite[Theorem 2.7]{K-S} tells us that there is an $\epsilon>0$, such
that any closed smooth Willmore sphere with Willmore functional $4\pi+\epsilon$
is a round sphere.  Therefore, at least one of $\widetilde{y}_1,\dots, \widetilde{y}_5$
is  a nonremovable singular point. However, a Willmore sphere cannot have only one singular point, by Lemma 4.2 in \cite{K-S2} (which is true in $\R^n$), therefore $\widetilde{f}$ has at least 2 singular points.
\endproof

\subsection{Minimizing Willmore functional subject to area constraint}

In this subsection, $N$ stands for a compact closed submanifold of $\R^n$ with induced metric.
We say $f\in W^{2,2}_{conf}(\Sigma,h,N)$ if $f\in W^{2,2}_{conf}
(\Sigma,h,\R^n)$ and $f(\Sigma)\subset N$. For $f\in W^{2,2}_{conf}(\Sigma,h,N)$, we define
$$
W(f)=W(f,\Sigma,N)=\frac{1}{4}\int_{\Sigma}|H_{f,\Sigma,N}|^2d\mu_f.
$$

First, we consider the case of genus zero. Set
$$
\beta_0(N,a)=\inf\{W(f):\mu(f)=a,\s f\in W^{2,2}_{conf}(S^2,N)\}.
$$

\begin{pro}  We have
$$\lim_{a\rightarrow 0}\beta_0(N,a)=4\pi.$$
Moreover,  when $a$ is sufficiently small, there is
an embedding $f\in W^{2,2}_{conf}(S^2,N)$, such that
$$\mu(f)=a,\s and\s W(f)=\beta_0(N,a).$$
\end{pro}

\proof First, we show that
\begin{equation}\label{beta0}
\limsup_{a\rightarrow 0}\beta_0(N,a)\leq 4\pi.
\end{equation}
Take a point $p\in N$ and a normal
coordinate neighborhood $U$ around $p$.
Let $$S_r=\{(x^1,x^2,x^3,0,\dots,0)\in T_pN:(x^1)^2+(x^2)^2+(x^3)^2=r^2\}.$$
It is easy to check that
$$\lim_{r\rightarrow 0}W(\exp_p(S_r),N)=4\pi.$$
For any $a$ which is sufficiently small, we can
find $r=r(a)$ such that $\mu(\exp_p(S_r))=a$
 and $r\rightarrow 0$ as $a\rightarrow 0$. Then \eqref{beta0} follows from
$\beta_0(N,a)\leq W(\exp_p(S_r))$.

Next, we prove that $\beta_0(N,a)$ can be attained
by an embedded 2-sphere.
Let $f_k\in W^{2,2}_{conf}(S^2,N)$ be a minimizing sequence of $\beta_0(N,a)$.
By Lemma \ref{H} and \eqref{beta0}, when $a$ is sufficiently small
and $k$ is sufficiently large
$$
W(f_k,S^2,\R^{n})\leq W(f_k,S^2,N)+C\mu(f_k)<4\pi + \epsilon(a,k)+Ca
$$
where $\epsilon(a,k)\to 0$ as $a\to 0$ and $k\to\infty$.
By Theorem \ref{main}, $\{f_k\}$ has a  limit $f_0$, which
is a branched conformal immersion from a stratified sphere
$S$ into $N$ with
$$\mu(f_0)=a\s \hbox{and}\s W(f_0)\leq \beta_0(N,a).$$
Then by (\ref{density1}), for any $p\in S$ it holds
$$
\theta^2(f_0(p))<2.
$$
Thus $S$ is a 2-sphere and $f_0$
has no branch points and no self-intersection points.
Hence $f_0$ is an embedding. Therefore $f_0$ is a minimizer for $\beta_0(N,a)$:
$$
W(f_0)=\beta_0(N,a).
$$

Finally, we prove
$$\varliminf_{a\rightarrow 0}\beta_0(N,a)\geq 4\pi.$$
By Lemma \ref{H},
$$W(f_0,S^2,\R^n)\leq W(f_0,S^2,N)+Ca.$$
It is well-known that $W(f_0,S^2,\R^n)\geq 4\pi$,
which completes the proof.
\endproof

We now consider the case of genus larger than 0.
Recall a result of Schoen-Yau \cite{S-Y} and Sacks-Uhlenbeck \cite{S-U2}:
If $\varphi:\Sigma\to N$ induces an injection from the fundamental groups to $\Sigma$ and $N$, then there is a branched minimal immersion $f: \Sigma\rightarrow N$ so that $f$ induces the same map between fundamental groups as $\varphi$ and $f$ has least area among all such maps. If $\pi_2(N)=0$ then $f$ is minimizing in its homotopy class.
We denote the area of the branched minimal immersion $f_\varphi$ by $a_\varphi$.

Let $g>0$ be the genus of the closed Riemann surface $\Sigma$ and $\phi:\Sigma\to N$ be a continuous map. Define
$$
\beta_g(N,a,\phi)=\inf\left\{W(f):f\in \widetilde{W}^{2,2}(\Sigma,N), \,\,\mu(f)=a, \,\,f\sim \phi\right\},
$$
where $f\sim\phi$ means that $f$ is homotopic
to $\phi$.

\begin{pro}  Let $\Sigma$ be a closed Riemann surface with genus $g>0$ and $N$ be a compact Riemannian  manifold with $\pi_2(N)=0$. Let
$\varphi:\Sigma\to N$ be a map which induces an injective $\varphi_\#:\pi_1(\Sigma)\to\pi_1(N)$.
 Then we can find an $\delta>0$, such that for any
$a\in [a_\varphi,a_\varphi+\delta)$,
there is a branched conformal immersion $f_0$ of a smooth Riemann surface $(\Sigma,h)$ of genus $g$ in $\R^n$, such that
$\mu(f_0)=a$ and $W(f_0)=\beta_g(N,a,\varphi)$ and $f_0$ is homotopic to $\varphi$.
Moreover, when $\dim N=3$, we can choose $\delta$
to be small such that $f_0$ is an immersion.
\end{pro}

\proof The proof will be divided into several steps.

{\bf Step 1.} We prove that
$\lim_{a\rightarrow a_0}\beta_g(N,a,\varphi)=0.$

Let $F\in C^\infty(\Sigma\times[0,1],\R^n)$,
such that $F(\cdot,t)$ is an immersion for each $t$  and
$$
F(\cdot,0)=f_\varphi,\s \mu(F(\cdot,1))\geq a_\varphi.
$$
As $F(\cdot,t)\sim\varphi$ and $f_\varphi$ is a minimal surface,
$$
\lim_{a\rightarrow a_\varphi}\beta_p(N,a,\varphi)\leq\lim_{t\rightarrow 0}W(F(\cdot,t))=W(f_\varphi)=0.
$$

{\bf Step 2.} Smooth convergence of conformal structures.

We take a minimizing sequence $\{f_k\}$ of
$\beta_g(N,a,f)$. Recall that $f_k$
are $W^{2,2}$  branched conformal immersions from $(\Sigma,h_k)$
into $\R^n$, where $h_k$ are  the smooth metrics
with curvature 0 or $-1$.
Because $\pi_2(N)=0$ and $f_k\sim\varphi$ for each $k$, $f_k$ induces the same injective action on the fundamental groups as $\varphi$ does; hence the conformal structures of $h_k$ stay in a compact set of the moduli space for both the cases $g>1$ and $g=1$, therefore, after passing to a subsequence if necessary, $\Sigma_k=(\Sigma,h_k)$ converges to a Riemann surface $(\Sigma,h_0)$ in $\mathcal{M}_g$ (cf. \cite{S-Y}). The results in \cite{S-Y} applies as $f_k$ belong to $W^{1,2}\cap C^0$.




{\bf Step 3.} We prove that $\{f_k\}$ has no bubbles, i.e. the limit $f_0$ is a map defined on $\Sigma$.

By Remark \ref{moduli}, $f_0$ is defined on $\Sigma_\infty=\Sigma_0\cup S_1\cup S_2\cdots\cup S_m$, where $S_i$ are all 2-spheres and $\Sigma_0$ is a smooth surface of genus $g$. We prove $m=0$. Assume $m\geq1$.
By Theorem \ref{main}, $\mu(f_0)=a$ and $W(f_0)\leq \beta_p(N,a,\varphi)$. Further, $f_k(\Sigma)$ converge to $f_0(\Sigma_\infty)$ in Hausdorff distance and $f_0|_{S_j}$ is homotopic to a constant map for each $j=1,\dots,m$
as $\pi_2(N)=0$. We conclude that $f_0|_{\Sigma_0}$ is homotopic to $\varphi$. Consequently, $\mu(f_0(\Sigma_0))\geq a_\varphi$. Then we get
$$
\mu(f_0,S_i)\leq\mu(\Sigma)-\mu(\Sigma_0)\leq a-a_\varphi \s\hbox{and}\s W(f_0,S_i,N)\leq \beta_g(N,a,\varphi).
$$
By Lemma \ref{H} and Step 1,
$$
W(f_0,S_i,\R^n)\leq C(a-a_\varphi)+\beta_g(N,a,\varphi)\to 0\s\hbox{as}\s a\to a_\varphi.
$$
This, however, contradicts Proposition \ref{gap} when the Willmore functional of $S_i$ goes below the gap constant.

{\bf Step 4.} We consider the case of $\dim N=3$.

We will use the result that there are no branch points for minimal surfaces
\cite{G,O} to prove that $f_0$ has no branch points when
$\delta$ is sufficiently small.

If the claimed result is not true, then  there is a sequence of numbers $a_k>a_\varphi$ with $a_k\rightarrow a_\varphi$ and a sequence of $W^{2,2}$ branched conformal immersions ${f}_{0,k}$ of $(\Sigma,h_k)$ in $N$ with
$\mu(f_{0,k})=a_k$, $W(f_{0,k},\Sigma,N)=\beta_g(N,a_k,\varphi)$ by the first part of the proposition, and each $f_{0,k}$ has at least a branch point $p_k$. By Step 1, $W(f_{0,k},\Sigma,N)\rightarrow 0$.

As in Step 2,  $(\Sigma,h_k)$ converge to a
smooth surface $(\Sigma,h_0)$ in $\mathcal{M}_g$.
For simplicity, we will still denote ${f}_{0,k}\circ\psi_k$ (see Remark \ref{moduli})
by ${f}_{0,k}$ which is a branched conformal immersion
from $(\Sigma,\psi^*_k(h_k))$ into $\R^n$.
By Theorem \ref{main1}, we may set
${f}_{0,0}$ to be the limit of ${f}_{0,k}$
with $\mu({f}_{0,0})=a$ and $W({f}_{0,0})=0$.
Arguing as in Step 3, $\{{f}_{0,k}\}$
has no bubbles, and ${f}_{0,0}\in W^{2,2}_{b,c}(\Sigma,h_0,\R^n)$ for some smooth $h_0$. Moreover, ${f}_{0,0}$ is a minimal surface in $N$. By the result of Gulliver and Osserman, ${f}_{0,0}$ is
a smooth immersion of $\Sigma$ in $N$.

Since $p_k$ is a branch point, by Theorem \ref{removal}, the area density
$$
\theta^2_{f_{0,k}(p_k)}(f_{0,k}(U))\geq 2,
$$
where $U$ is a neighborhood of $p_k\rightarrow p$ in $\Sigma$ for sufficiently large $k$.
As $f_{0,0}$ is immersive, we can take $U$ small so that
${f}_{0,0}$ is an embedding on $U$ and $\mu({f}_{0,0}(U))<\epsilon'$.
Further, by the monotonicity formula for minimal surfaces, for small $r$ and geodesic balls $B^N_r(f_{0,0}(p))$ in $N$, it holds
$$
\mu({f}_{0,0}(U)\cap B^N_r({f}_{0,0}(p)))\leq (1+\epsilon')\pi r^2.
$$
From the expansion of metric in normal coordinates, for small $r$ and the Euclidean ball
$B_{r}(f_{0,0}(p))$ in $\R^n$ we have
$$
\mu(f_{0,0}(U)\cap B_r(f_{0,0}(p)))   \leq \mu( f_{0,0}(U)\cap B^N_{r+cr^2}(f_{0,0}(p)))\leq (1+\epsilon')\pi r^2+O(r^3)
$$
where $c$ depends on $N$.

In light of  Lemma \ref{H}, $W({f}_{0,k},U,\R^n)<\epsilon_0^2$
if we choose $\epsilon'$ to be very small and $k$ large enough. Then
$\{{f}_{0,k}\}$ has no blowup points in $U$ by the $\epsilon$-regularity.
Then we have
$$
\mu({f}_{0,k}(U)\cap B_r({f}_{0,k}(p_k)))\rightarrow
\mu({f}_{0,0}(U)\cap B_r({f}_{0,0}(p)))\s\s\hbox{as $k\to\infty$}.
$$
By Lemma \ref{H},
$$
W({f}_{0,k},U,\R^n)\leq C\epsilon'+W(f_{0,k},U,N).
$$
Then by \eqref{density2},
\begin{eqnarray*}
\theta^2_{f_{0,k}(p)}(f_{0,k}(U))&\leq& \frac{\mu(f_{0,k}(U)\cap B_r(f_{0,k}(p_k)))}{\pi r^2}+
W(f_{0,k},U,\R^n)+CW(f_{0,k},U,\R^n)^{\frac{1}{2}}.
\end{eqnarray*}
Hence,
\begin{eqnarray*}
2&\leq& \lim_{U\to p}\lim_{k\to\infty}\left(
\frac{\mu(f_{0,k}(U)\cap B_r(f_{0,k}(p_k)))}{\pi r^2}+
W(f_{0,k},U,\R^n)+CW(f_{0,k},U,\R^n)^{\frac{1}{2}}\right)\\
&\leq& 1+\epsilon'.
\end{eqnarray*}
This is impossible for $\epsilon'$ small.
\endproof

\subsection{Minimizing Willmore functional of surfaces with a Douglas type condition}
In this subsection, we consider a sufficient condition of Douglas type as in the minimal surface theory for
existence of minimizers of the Willmore functional.

First, we assume $N$ to be a compact Riemannian manifold with negative sectional curvatures.
In negatively curved $N$, surface area is bounded by the Willmore functional and the genus of the surface.

\begin{lem}\label{negative.curvature} Let $N$ be a compact Riemannian manifold
with $K\leq -c<0$. Then for any $f\in \widetilde{W}^{2,2}(\Sigma,N)$,
$$\mu(f)\leq c^{-1}\left(W(f,\Sigma,N)-2\pi\chi(\Sigma)\right).$$
Especially, when $g(\Sigma)=0$ or $1$,
$$\mu(f)\leq c^{-1}W(f,\Sigma,N).$$
\end{lem}

\proof

From the Gauss equation:
$$
R^\Sigma(X,Y,X,Y)=R^N(X,Y,X,Y)+\lan A(X,X),A(Y,Y)\ran-\lan A(X,Y),A(X,Y)\ran.
$$
we have
$$
K_\Sigma\leq K_{f_*(T\Sigma)}+\frac{1}{4}|H_{f,\Sigma,N}|^2.
$$
Then from the generalized Gauss-Bonnet formula - Lemma \ref{GB}, we have
$$
2\pi\chi(\Sigma)+2\pi b\leq -c\mu_{f}(\Sigma)+W(f,\Sigma,N)
$$
where $b$ is the number of branch points, in turn
$$
c\mu_f(\Sigma)\leq W(f,\Sigma,N)-2\pi\chi(\Sigma).
$$
When $g(\Sigma)\leq1$ the Euler number $\chi(\Sigma)$ is nonnegative, in this case
$$
c\mu_{f}(\Sigma)\leq W(f,\Sigma,N).
$$
Dividing by $c$ yields the desired area bounds. \endproof

Recall that any connected stratified surface $\Sigma$ can be written as union of finitely many connected 2-dimensional components: $\Sigma=\bigcup_i\Sigma_i$. Denote the genus of $\Sigma$ and $\Sigma_i$ by
$g(\Sigma)$ and $ g(\Sigma_i)$, accordingly.
We introduce a subset $S(g)$ of all stratified surfaces as follows.
\begin{enumerate}
\item If $g>0$,
$
S(g) =\left \{\Sigma: \hbox{$\Sigma=\bigcup_i \Sigma_i$  with  $g(\Sigma_i) < g$ for all $i$}\right\}.
$
\item If $g=0$,
$
S(0)=\left\{ \Sigma: \hbox{$\Sigma=\bigcup_i \Sigma_i$  with $g(\Sigma)=0$ and $i\geq 2$}\right\}.
$
\end{enumerate}
Note that any $\Sigma\in S(g)$ with $g(\Sigma)=g$ must be singular, in the sense that it has more than one components. Especially, $S(g)\cap{\mathcal M}_g=\emptyset$. However, when $g\geq 1$, $S(g)$ contains smooth surfaces of genus $\leq g-1$.



Define
\begin{eqnarray*}
\alpha^*(g)&=&\inf \{ W(f,\Sigma,\R^n): f\in W^{2,2}_{b,c}(\Sigma,\R^n), f(\Sigma)\subset N,\Sigma\in S(g)\} \\
\alpha(g)&=& \inf \{W(f,\Sigma,\R^n): f\in W^{2,2}_{b,c}(\Sigma,\R^n), f(\Sigma)\subset N,\Sigma\in{\mathcal M}_g\}.
\end{eqnarray*}


We now state a sufficient condition, similar to the Douglas condition for minimal surfaces, for existence of minimizers for the Willmore functional.

\begin{pro}
Let $N$ be a compact Riemannian manifold with negative sectional curvatures. If $0<\alpha(g)<\alpha^*(g)$, then there is
a $W^{2,2}$ branched conformal immersion $f$ from a smooth surface of genus $g$ into $N$ which minimizes the Willmore functional among all such maps.
\end{pro}

\proof Let $f_k:(\Sigma,h_k)\rightarrow N\hookrightarrow\R^n$ be a minimizing sequence of $\alpha(g)$.
By
Lemma \ref{negative.curvature},
the areas $\mu(f_k(\Sigma))$ are uniformly bounded as well by Lemma \ref{negative.curvature}. Since $\alpha(g)$ is positive, $f_k$ cannot converge to a point. Then from Theorem \ref{main},  there exists a subsequence of $\{f_k\}$, still denoted by $\{f_k\}$, a limit map $f_0\in W^{2,2}_{b,c}(\Sigma_\infty,\R^n)$ from a stratified Riemann surface $\Sigma_\infty$ with $g(\Sigma_\infty)\leq g$ into $N\hookrightarrow\R^n$, and
$$
W(f_0,\Sigma_\infty,\R^n)\leq \lim_{k\to\infty}W(f_k,\Sigma,\R^n) = \alpha(g).
$$

We write $\Sigma_\infty=\bigcup_{i=1}^m \Sigma_i$. If $g(\Sigma_\infty)=g$, we consider two cases.
Case 1: $g(\Sigma_i)=g$ for some $i=1,...,m$. In this case,
$$
 W(f_0|_{\Sigma_1},\Sigma_1,\R^n)\leq W(f_0,\Sigma_\infty,\R^n)=\alpha(g).
$$
So $f_0(\Sigma_i)$ is a smooth genus $g$ surface attains $\alpha(g)$. Case 2:  $g(\Sigma_i)<g$ for all $i=1,...,m$.
Thus $\Sigma_\infty\in S(g)$, and in turn
$$
\alpha^*(g)\leq W(f_0,\Sigma_\infty,\R^n)\leq \alpha(g)<\alpha^*(g).
$$
This contradiction rules out Case 2.
If $g(\Sigma_\infty)<g$ then $\Sigma_g\in S(g)$. Therefore
$$
\alpha^*(g)\leq W(f_0,\Sigma_\infty,\R^n)\leq\alpha(g)<\alpha^*(g)
$$
and this is impossible.
\endproof

Instead of the curvature assumption on $N$, we set, for $0<a<\infty$,
\begin{eqnarray*}
\gamma^*(g,a)&=&\inf \{ W(f,\Sigma,\R^n): f\in W^{2,2}_{b,c}(\Sigma,\R^n), f(\Sigma)\subset N,\Sigma\in S(g), \mu(f(\Sigma))\leq a\} \\
\gamma(g,a)&=& \inf \{W(f,\Sigma,\R^n): f\in W^{2,2}_{b,c}(\Sigma,\R^n), f(\Sigma)\subset N,\Sigma\in{\mathcal M}_g,
\mu(f(\Sigma))\leq a\}.
\end{eqnarray*}
Since there is no loss in measures in the limit process, as asserted in Theorem \ref{main}, the same proof above allows us to conclude

\begin{pro}
Let $N$ be a compact Riemannian manifold. If $0<\gamma(g,a)<\gamma^*(g,a)$, then there is
a $W^{2,2}$ branched conformal immersion $f$ from a smooth surface of genus $g$ into $N$ which minimizes the Willmore functional among all such maps.
\end{pro}

\section{Appendix}

Wente's inequality \cite{B,Ge,To} states that if $u\in W^{1,2}_0(D)$
solves the equation
$$-\Delta u=\nabla a\,\nabla^\bot b,$$
then we have
\begin{equation}\label{W1}
\|u\|_{L^\infty(D)}\leq\frac{1}{2\pi}\|\nabla a\|_{L^2(D)}
\|\nabla b\|_{L^2(D)}.
\end{equation}
and
\begin{equation}\label{W2}
\|\nabla u\|_{L^2(D)}
\leq \frac{3}{16\pi}\|\nabla a\|_{L^2(D)}\|\nabla b\|_{ L^2(D)}.
\end{equation}

\begin{lem}\label{u-conti}
Let $u\in W^{1,2}_0(D)$ be the unique solution to the equation
$$-\Delta u=\nabla a\nabla^\bot b,$$
where $a,b\in W^{1,2}(D)$. Then $u\in C^0(D)$.
\end{lem}

\proof

Let $a_k\in C^\infty(\bar{D})$ with $a_k\rightarrow a$ in $W^{1,2}(D)$.
There exist solutions $u_k\in W^{2,2}(D)\cap C^{0,\alpha}(D)$ to the Dirichlet problem
\begin{eqnarray*}
-\Delta u_k&=&\nabla a_k\nabla^\bot b,\,\,\,\,\s\s \mbox{in $D$ }\\
 u_k&=&0, \,\,\,\,\,\,\s\s\s\s\s\s\mbox{on $\partial D$.}
\end{eqnarray*}
We have
$$-\Delta(u_k-u_m)=\nabla (a_k-a_m)\nabla^\bot b.$$
Then by \eqref{W1}
\begin{eqnarray*}
\|u_k-u_m\|_{C^0(D)}&=&\|u_k-u_m\|_{L^\infty(D)}\\
&\leq& C\,\|a_k-a_m\|_{L^2(D)}\|b\|_{L^2(D)}.
\end{eqnarray*}
Hence $u_k$ converge in $C^0(D)$ to a continuous function $u_0$ which vanished on $\partial D$.
By \eqref{W2},
we may assume $u_k$ converges to $u_0$ in $W^{1,2}(D)$ as well.
For any smooth $\varphi$,
$$\int_D\varphi\nabla a\nabla^\bot b=\lim_{k\rightarrow+\infty}
\int_D\varphi\nabla a_k\nabla^\bot b
=\lim_{k\rightarrow+\infty}
\int_D\nabla u_k\nabla \varphi=\int_D\nabla u_0\nabla \varphi.$$
Hence $-\Delta u_0=\nabla a\nabla^\bot b$, and it follows from the uniqueness of solution to the Dirichlet problem
that  $u=u_0$, so $u\in C^0(D)$.
\endproof



Now, we are ready to prove:
\begin{pro}  If $f\in W^{2,2}_{conf}(D,\R^n)$ with
$df\otimes df=e^{2u}g_{euclid}$ and $\| u \|_{L^\infty}<+\infty$, then $u\in C^0(D)$.
\end{pro}

\proof
In a complex coordinate $z=x+iy$ on $D$, since $f$ is conformal and the induced metric is $e^{2u}(dx^2+dy^2)$, we have
$$
f_x\cdot f_x = f_y\cdot f_y =e^{2u}
$$
and we can take $a =e^{-u} \,f_x, b=e^{-u}\, f_y$ as a local orthonormal frame for the tangent bundle of $f(D)$. Straight computation shows
$$
a_x = e^{-u}f_{xx}+(e^{-u})_x f_x=e^{-u}f_{xx}-e^{-3u}(f_{xx}\cdot f_x)f_x
$$
which leads to
$$
a_x\cdot a_x \leq C e^{-2u} f_{xx}\cdot f_{xx}\leq C e^{2\|u\|_{L^\infty}}\left| f_{xx}\right|^2
$$
and similarly
$$
a_y\cdot a_y \leq  C e^{2\|u\|_{L^\infty}}\left| f_{yy}\right|^2.
$$
Therefore, $a\in W^{1,2}(D)$ and similarly $b\in W^{1,2}(D)$. On the other hand, we can check
$$
K_fe^{2u}=\nabla a\cdot\nabla^\perp b,\s -\Delta u=K_fe^{2u}.
$$
Let $v$ be the harmonic function on $D$ which agrees with $u$ on $\partial D$. Then
\begin{eqnarray*}
\Delta (u-v)&=& \nabla a \cdot \nabla^\perp b\,\,\,\,\,\s\hbox{in $D$}\\
u-v&=&0\,\,\,\,\,\,\s\s\s\s\s\s\hbox{on $\partial D$}
\end{eqnarray*}
and Lemma \ref{u-conti} implies $u-v\in C^0(D)$, hence $u\in C^0(D)$. \endproof


\vspace{2ex}

\begin{tabular}{ll}
{\sc Jingyi Chen }&{\sc Yuxiang Li}\\
{\sc Department of Mathematics}&{\sc Department of Mathematical Sciences}\\
{\sc University of British Columbia}&{\sc Tsinghua University}\\
{\sc Vancouver,B.C., V6T1Z2, Canada}& {\sc Beijing 100084, P.R. China}\\
{\tt jychen@math.ubc.ca} & {\tt yxli@math.tsinghua.edu.cn}
\end{tabular}


\begin{thebibliography}{2}

\bibitem{B}S. Baraket:
Estimations of the best constant involving the $L^\infty$
norm in Wente's inequality,
{\em Ann. Fac. Sci. Toulouse Math. (6)} {\bf 5} (1996), 373¨C385.



\bibitem{C-T} J. Chen and G. Tian: Compactification of moduli space
of harmonic mappings, \emph{ Comment. Math. Helv.}  {\bf 74}  (1999), 201-237.

\bibitem{C-L-W} L. Chen, Y. Li, Y. Wang: The Refined Analysis on the Convergence Behavior of Harmonic Map Sequence from Cylinders, \emph{ J. Geom. Anal.}, to appear.




\bibitem{D-K}D. DeTurck and J. Kazdan:
Some regularity theorems in Riemannian geometry,
\emph{Ann. Sci. \'Ecole Norm. Sup. (4)} {\bf 14}
 (1981),  249-260.

\bibitem{D-T} W. Ding and G. Tian: Energy identity for a class of approximate
harmonic maps from surfaces, \emph{Comm. Anal. Geom.} {\bf 3} (1995),
543-554.



\bibitem{E-T} J. Eschenburg and R. Tribuzy:
Branch points of conformal mappings of surfaces, \emph{Math. Ann.} {\bf 279} (1988), 621-633.

\bibitem{Ge}Y. Ge:
Estimations of the best constant involving the
$L^2$ norm in Wente's inequality and compact
$H$-surfaces in Euclidean space,
{\em ESAIM Control Optim. Calc. Var.} {\bf 3} (1998),
263--300.

\bibitem{G} R. Gulliver: Regularity of minimizing surfaces of prescribed mean curvature,
 \emph{ Ann. of Math. (2)} {\bf 97} (1973), 275-305.

\bibitem{Halpern} N. Halpern: A proof of the collar lemma,
\emph{Bull. London Math. Soc.}  {\bf 13}  (1981),  141-144.

\bibitem{He} F. H\'elein: Harmonic maps,
conservation laws and moving frames.
Translated from the 1996 French original.
With a foreword by James Eells. Second edition.
Cambridge Tracts in Mathematics, 150.
Cambridge University Press, Cambridge, 2002.



\bibitem{Hum} C. Hummel: Gromov's compactness
theorem for pseudo-holomorphic curves,
 \emph{Progress in Mathematics} {\bf 151}, Birkh\"auser Verlag, Basel (1997).

\bibitem{Keen} L. Keen: Collars on Riemann surfaces,
Discontinuous groups and Riemann surfaces (Proc. Conf., Univ. Maryland, College Park, Md., 1973),  263268. Ann. of Math. Studies, No. 79, Princeton Univ. Press, Princeton, N.J., 1974.

\bibitem{K-L} E. Kuwert and Y. Li: $W^{2,2}$-conformal immersions of a closed Riemann
surface into $\R^n$,  \emph{arXiv:1007.3967}.

\bibitem{K-S}
E. Kuwert and R. Sch\"atzle:
The Willmore flow with small initial energy,
\emph{ J. Differential Geom.}, {\bf 57} (2001), 409-441.



\bibitem{K-S2} E. Kuwert and R. Sch\"atzle: Removability
of point singularities of Willmore surfaces, \emph{ Ann. of Math.}
{\bf 160} (2004), 315-357.

\bibitem{La} T. Lamm: Energy identity for approximations of harmonic maps from surfaces,
  \emph{Trans. Amer. Math. Soc.}  {\bf 362}  (2010),   4077-4097.


\bibitem{La-M} T. Lamm and J. Metzger: Small surfaces of Willmore type in Riemannian manifolds,
 \emph{Int. Math. Res. Not.} (2010), no.  3786-3813.

\bibitem{L-W} Y. Li and Y. Wang: A weak energy identity and the length of necks for a Sacks-Uhlenbeck -harmonic map sequence, \emph{ Adv. Math.} {\bf 225} (2010), 1134-1184 .





\bibitem{Lin-W} F. Lin and  C. Wang:  Energy identity of harmonic
map flows from surfaces at finite singular time,
 \emph{Calc. Var. Partial Differential Equations}, {\bf 6}  (1998),   369-380.

\bibitem{M} J.P. Matelski: A compactness theorem for Fuchsian groups of the second kind,
\emph{Duke Math. J.} {\bf 43} (1976),  829-840.


\bibitem{M-S} S. M\"uller and V. \v{S}ver\'ak: On surfaces
of finite total curvature, \emph{ J. Differential Geom.}
{\bf 42} (1995),
229-258.

\bibitem{O} R. Osserman:
A proof of the regularity everywhere of the classical solution to Plateau's problem,
\emph{ Ann. of Math. (2)} {\bf 91} (1970) 550-569.

\bibitem{P} T. H. Parker: Bubble tree convergence
for harmonic maps,
\emph{ J. Differential Geom.} {\bf 44} (1996), 595-633.

\bibitem{Q} J. Qing:  On singularities of the heat flow for harmonic maps from surfaces into
spheres, \emph{Comm. Anal. Geom.}, {\bf 3} (1995), 297-315.

\bibitem{Q-T} J. Qing and  G. Tian: Bubbling of the heat flow for harmonic maps from surfaces,
\emph{Comm. Pure. Apple. Math.}, {\bf 50} (1997), 295-310.

\bibitem{Ra} B. Randol: Cylinders in Riemann surfaces,
\emph{Comment. Math. Helv.} {\bf 54} (1979),  1-5.

\bibitem{R}T. Rivi\'ere: Analysis aspects of Willmore
surfaces, {\em Invent. Math.} {\bf 174} (2008), 1-45.


\bibitem{S-U1} J. Sacks and K.  Uhlenbeck: The existence of minimal immersions of
2-spheres, \emph{ Ann. of Math.}, {\bf 113} (1981), 1-24.

\bibitem{S-U2} J. Sacks and K. Uhlenbeck: Minimal immersions of closed Riemann surfaces, \emph{Trans. Amer. Math. Soc.} {\bf 271} (1982), no. 2, 639-652.

\bibitem{S-Y} R. Schoen and S.T. Yau,
Existence of incompressible minimal surfaces and the topology of three-dimensional manifolds with nonnegative scalar curvature,
\emph{Ann. of Math. (2)} {\bf 110} (1979), no. 1, 127-142.

\bibitem{S} L. Simon: Existence of surfaces minimizing
the Willmore functional, \emph{Comm. Anal. Geom.}
{\bf 1} (1993),
281-326.

\bibitem{Sh} R. Sch\"atzle: The Willmore boundary problem,
\emph{Cal. Var. P.D.E.}, {\bf 37} (2010),  275-302.


\bibitem{To} P. Topping:  The optimal constant in Wente's $L^\infty$
 estimate, {\em Comment. Math. Helv.} {\bf 72} (1997), 316¨C328.


\bibitem{Z} M. Zhu: Harmonic maps from degenerating Riemann surfaces,
\emph{Math. Z.}  {\bf 264}  (2010),  no. 1, 63-85.

\end{thebibliography}
\end{document}